\documentclass[12pt]{amsart}
\usepackage{enumerate, amssymb}
\usepackage[matrix,arrow,curve]{xy}
\usepackage{mathrsfs}

\input{diagrams.tex}
\def\pdfsyncstart{}
\def\pdfsyncstop{}
\theoremstyle{plain}

\newtheorem{thm}{Theorem}[section]
\newtheorem{cor}[thm]{Corollary}
\newtheorem{lem}[thm]{Lemma}
\newtheorem{prop}[thm]{Proposition}

\theoremstyle{definition}
\newtheorem{defi}[thm]{Definition}
\newtheorem{defis}[thm]{Definitions}
\newtheorem{conj}[thm]{Problem}
\newtheorem{conv}[thm]{Convention}
\newtheorem{nota}[thm]{Notation}
\newtheorem{rem}[thm]{Remark}
\newtheorem{rems}[thm]{Remarks}
\newtheorem{exa}[thm]{Example}
\newtheorem{exas}[thm]{Examples}
\newtheorem{sit}[thm]{}

\newcommand{\brem}{\begin{rem}}
\newcommand{\brems}{\begin{rems}}
\newcommand{\erem}{\end{rem}}
\newcommand{\erems}{\end{rems}}
\newcommand{\bexa}{\begin{exa}}
\newcommand{\bexas}{\begin{exas}}
\newcommand{\eexa}{\end{exa}}
\newcommand{\eexas}{\end{exas}}
\newcommand{\bdefi}{\begin{defi}}
\newcommand{\edefi}{\end{defi}}
\newcommand{\bdefis}{\begin{defis}}
\newcommand{\edefis}{\end{defis}}
\newcommand{\bcor}{\begin{cor}}
\newcommand{\ecor}{\end{cor}}
\newcommand{\blem}{\begin{lem}}
\newcommand{\elem}{\end{lem}}
\newcommand{\bconv}{\begin{conv}}
\newcommand{\econv}{\end{conv}}
\newcommand{\bconj}{\begin{conj}}
\newcommand{\econj}{\end{conj}}
\newcommand{\bprop}{\begin{prop}}
\newcommand{\eprop}{\end{prop}}
\newcommand{\bthm}{\begin{thm}}
\newcommand{\ethm}{\end{thm}}
\newcommand{\bnota}{\begin{nota}}
\newcommand{\enota}{\end{nota}}
\newcommand{\bsit}{\begin{sit}}
\newcommand{\esit}{\end{sit}}
\newcommand{\be}{\begin{equation}}
\newcommand{\ee}{\end{equation}}
\newcommand{\bproof}{\begin{proof}}
\newcommand{\eproof}{\end{proof}}
\def\ba{\begin{array}}
\def\ea{\end{array}}
\def\bea{\begin{eqnarray}}
\def\eea{\end{eqnarray}}

\def\bnum{\begin{enumerate}}
\def\enum{\end{enumerate}}

\newcommand{\no}{\noindent}

\newcommand{\Der}{\operatorname{Der}}
\newcommand{\Rad}{\operatorname{Rad}}

\newcommand{\Aut}{{\operatorname{Aut}}}
\newcommand{\Lin}{{\operatorname{Lin}}}
\newcommand{\Bir}{{\operatorname{Bir}}}
\newcommand{\Alb}{{\operatorname{Alb}}}

\newcommand{\Frac}{\operatorname{Frac}}
\newcommand{\spec}{\operatorname{spec}}
\newcommand{\supp}{\operatorname{supp}}
\newcommand{\mult}{\operatorname{mult}}
\newcommand{\divis}{\operatorname{div}}
\newcommand{\SL}{{\operatorname{SL}}}
\newcommand{\Bl}{\operatorname{Bl}}
\newcommand{\affcone}{\operatorname{AffCone}}
\newcommand{\Pic}{\operatorname{Pic}}
\newcommand{\Exc}{\operatorname{Exc}}

\newcommand{\GL}{\operatorname{GL}}
\newcommand{\ML}{\operatorname{ML}}
\newcommand{\PGL}{\operatorname{PGL}}
\newcommand{\PSO}{\operatorname{PSO}}

\newcommand{\Spec}{\operatorname{Spec}}
\newcommand{\Div}{\operatorname{Div}}
\newcommand{\pr}{\operatorname{pr}}
\newcommand{\Bs}{\operatorname{Bs}}

\newcommand{\A}{{\mathbb A}}
\newcommand{\PP}{{\mathbb P}}
\newcommand{\DD}{{\mathbb D}}

\newcommand{\C}{{\mathbb C}}
\newcommand{\Q}{{\mathbb Q}}
\newcommand{\Z}{{\mathbb Z}}
\newcommand{\N}{{\mathbb N}}
\newcommand{\T}{{\mathbb T}}
\newcommand{\V}{{\mathbb V}}
\newcommand{\F}{{\mathbb F}}

\newcommand{\G}{{\Gamma}}

\newcommand{\p}{{\partial}}

\newcommand{\bX}{{\bar X}}

\def\cH{{\mathcal H}}
\def\cL{{\mathcal L}}

\def\cO{{\mathcal O}}

\def\cH{{\mathscr H}}
\def\cL{{\mathscr L}}

\newcommand{\nlin}{\unitlength1mm\begin{picture}(0,9.25)
 \put(0,0.75){\line(0,1){8.5}}
 \end{picture}}

\newcommand{\vlin}[1]{\hspace{0.75mm}\unitlength1mm\begin{picture}
(#1,0)
 \put(0,0){\line(1,0){#1}}
 \end{picture}\hspace{0.75mm}\rule[-3mm]{0mm}{4mm}}

\def\llin{\vlin{11.5}}
\newcommand{\lin}{\vlin{8.5}}

\newcommand{\mybox}{\unitlength1mm\begin{picture}(0,1.5)
 \put(-0.75,-0.75){\line(0,1){1.5}}
 \put(-0.75,-0.75){\line(1,0){1.5}}
 \put(0.75,0.75){\line(0,-1){1.5}}
 \put(0.75,0.75){\line(-1,0){1.5}}
 \end{picture}}
\newcommand{\boxo}[1]{\unitlength1mm\begin{picture}(0,8)
 \put(0,0){\mybox}
 \put(0,3){\makebox(0,5)[b]{$#1$}}
 \end{picture}}

\newcommand{\cou}[2]{\unitlength1mm\begin{picture}(0,8)
 \put(0,0){\circle{1.5}}
 \put(0,3){\makebox(0,5)[b]{$#1$}}
 \put(0,-7){\makebox(0,4)[t]{$#2$}}
 \end{picture}
 \rule[-7mm]{0mm}{7mm}}

\newcommand{\crl}[2]{\unitlength1mm\begin{picture}(0,8)
 \put(0,0){\circle{1.5}}
 \put(-5,0){\makebox(0,5)[b]{$#1$}}
 \put(5,0){\makebox(0,5)[b]{$#2$}}
 \end{picture}
 \rule[-7mm]{0mm}{7mm}}

\newcommand{\cshiftup}[2]{\unitlength1mm\begin{picture}(0,9.25)
 \put(0,10){\crl{#1}{#2}}
 \end{picture}}

\newcommand{\boxrl}[2]{\unitlength1mm\begin{picture}(0,8)
 \put(0,0){\mybox}
 \put(-5,0){\makebox(0,5)[b]{$#1$}}
 \put(5,0){\makebox(0,5)[b]{$#2$}}
 \end{picture}
 \rule[-7mm]{0mm}{7mm}}

\newcommand{\boxshiftup}[2]{\unitlength1mm\begin{picture}(0,9.25)
 \put(0,10){\boxrl{#1}{#2}}
 \end{picture}}

\def\bdi{\pdfsyncstop\begin{diagram}}
\def\edi{\end{diagram}\pdfsyncstart}

\title[Group actions on affine cones]{Group actions on affine cones}

\author{Takashi Kishimoto}

\address{Department of Mathematics,
Faculty of Science, Saitama University, Saitama 338-8570, Japan}
\email{tkishimo@rimath.saitama-u.ac.jp}

\author{Yuri Prokhorov}

\address{Department
of Algebra, Faculty of Mathematics, Moscow State University,
Moscow 117234, Russia} \email{prokhoro@mech.math.msu.su}

\author{Mikhail Zaidenberg}

\address{Universit\'e
Grenoble I, Institut Fourier, UMR 5582 CNRS-UJF, BP 74, 38402 St.\
Martin d'H\`eres c\'edex, France} \email{zaidenbe@ujf-grenoble.fr}

\thanks{
The first author was supported by a Grant-in-Aid for Scientific
Research of JSPS No. 20740004. The second author was partially
supported by RFBR, No. \ 08-01-00395-a and Leading Scientific
Schools (grants No. \ NSh-1983.2008.1, NSh-1987.2008.1). This work
was done during a stay of the second and third authors at the Max
Planck Institut f\"ur Mathematik at Bonn and a stay of the first
and the second authors at the Institut Fourier, Grenoble. The
authors thank these institutions for hospitality.}

\begin{document}

\dedicatory{\footnotesize\it To Peter Russell on the occasion of
his 70th birthday}

\begin{abstract}
We address the following question:

\smallskip

\no {\it For which smooth projective varieties, the corresponding
affine cone admits an action of a connected algebraic group
different from the standard $\C^*$-action by scalar matrices and
its inverse action?}

\smallskip

We show in particular that the affine cones over anticanonically
embedded smooth del Pezzo surfaces of degree $\ge 4$ possess such
an action. Besides, we give some examples of rational Fano
threefolds which have this property. A question in \cite{FZ1}
whether this property holds also for smooth cubic surfaces, occurs
to be out of reach for our methods. Nevertheless, we provide a
general geometric criterion that could be helpful in this case as
well.
\end{abstract}

\maketitle

\bigskip

{\footnotesize \tableofcontents}

\section*{Introduction}
All varieties in this paper are defined over $\C$. By Corollary
1.13 in \cite{FZ1}, an isolated Cohen-Macaulay singularity $(X,x)$
of a normal quasiprojective variety $X$ is rational provided that
$X$ admits an effective action of the additive group $\C_+$, in
particular of a connected non-abelian algebraic group. In the
opposite direction, let us observe that, for instance, the
singularity at the origin of the affine Fermat cubic in $\A^4$
$$x_1^3+x_2^3+x_3^3+x_4^3=0\,$$ is rational.
The question was raised \cite[Question 2.22]{FZ1} whether it
also admits a non-diagonal action of a connected algebraic group,
in particular, a $\C_+$-action. So far, we do not know the answer.
However, we answer in affirmative a similar question for all del
Pezzo surfaces of degree $d\ge 4$.

\bthm\label{main1} Let $Y_d$ be a smooth del Pezzo surface of
degree $d$ anticanonically embedded into $\PP^d$, and let
$X_d=\affcone (Y_d)\subseteq\A^{d+1}$ be the affine cone over
$Y_d$. If $4\le d\le 9$ then $X_d$ admits a nontrivial
$\C_+$-action. Consequently, the automorphism group $\Aut\, (X_d)$
is infinite dimensional. Moreover, the Makar-Limanov invariant of
$X_d$ is trivial. \ethm

Recall \cite[10.1.1]{Dol1} that for $d\le 5$ the group $\Aut\,
(Y_d)$ is finite. By definition, the Makar-Limanov invariant of an
affine variety $X$ is the subring of $\cO(X)$ of common invariants
of all $\C_+$-actions on $X$. It is trivial when it consists of
the constants.

One of our main results (Theorem \ref{ext1}) provides a necessary
and sufficient condition for the existence of a nontrivial
$\C_+$-action on an affine cone. As a corollary, for affine cones
of dimension 3 we obtain the following geometric criterion.

\bthm\label{main2} Let $Y$ be a smooth projective rational surface
with a polarization $\varphi_{|H|}:Y\hookrightarrow \PP^n$, and
let $X=\affcone_H (Y)\subseteq\A^{n+1}$ be the affine cone over
$Y\subseteq\PP^n$. Then $X$ admits a nontrivial $\C_+$-action if
and only if $Y$ contains an $H$-polar cylinder i.e., a cylindrical
Zariski open set
$$U=Y\setminus \supp (D)\simeq Z\times \A^1,\,$$ where $Z$ is an
affine curve and $D\in |dH|$ is an effective divisor on $\PP^n$.
\ethm

Using this criterion, we show in Proposition \ref{ratsurfcon} that
for every smooth projective rational surface $Y$ there exists a
polarization $\varphi_{|H|}:Y\hookrightarrow \PP^n$ such that $Y$
contains an $H$-polar cylinder and so  the corresponding affine
cone possesses an effective action of $\C_+$. It would be
interesting to classify  in any dimension all pairs $(Y,H)$, where
$Y$ is a smooth projective variety and $H$ an ample divisor on
$Y$, such that the affine cone $X=\affcone_H (Y)$ admits an
effective $\C_+$-action. We recover this classification for
$\dim_\C (Y)=1$ and give some concrete examples in higher
dimensions, especially in dimensions $2$ and $3$.

A theorem due to Matsumura, Monsky and Andreotti (see \cite{MM},
or \cite[\S I.4]{GH}, or Section 1 below) claims that any
automorphism of a smooth hypersurface $Y$ in $\PP^n$ of degree
$d$, where $d,n\ge 3$ and $(d,n)\neq(4,3)$, is restriction of a
unique projective linear transformation, and $\Aut (Y)$ is a
finite group. In Corollary \ref{liext} we show that the
automorphism group $\Aut (X)$ of the affine cone $X$ over a
smooth, non-birationally ruled projective variety $Y$ is a linear
group, and actually a central extension of a finite group by
$\C^*$. Consequently, among the affine cones over smooth
projective surfaces in $\PP^3$, only those of degree $\le 3$ can
admit a nontrivial action of a connected algebraic group, and
their automorphism groups can be infinite-dimensional. Actually,
the 3-fold affine quadric cone possesses an effective linear
action of the additive group $\C_+^2$, see Example \ref{ex} below
or \cite{AS}, \cite{Sh}.

In Section 1 we give a short overview of the known results on the
automorphism groups. In Section 2 we collect generalities on
automorphisms of affine cones. Theorems \ref{main1} and
\ref{main2} are proven in Section 3. In Section 4 we summarize
some geometric facts that could be useful (in view of the
criterion of Theorem \ref{main2}) in order to answer Question 2.22
in \cite{FZ1} cited above. In the final section 5 we describe two
families of rational Fano threefolds such that the affine cones
over their anti-canonical embeddings possess effective
$\C_+$-actions.

{\bf Acknowledgements:}  We thank Alvaro Liendo, who participated
in discussions on early stages of this work, for his attention and
comments. Our thanks also to Dmitry Akhiezer, Michel Brion, and
Dimitri Timashev for providing useful information on homogeneous
varieties, to Dmitry Akhiezer and Ivan Arzhantsev for reading some
chapters and valuable remarks.

\section{Group actions on projective varieties}

In this section we recall some well known facts about the
automorphism groups of projective or quasiprojective varieties;
see e.g., \cite[\S I.4]{GH}, \cite[\S II.3]{LZ}. For an algebraic
variety $Y$, we let $\Aut (Y)$ denote the group of all biregular
automorphisms of $Y$ and $\Bir (Y)$ the group of all birational
transformations of $Y$ into itself. For a projective or an affine
embedding $Y\hookrightarrow \PP^n$ ($Y\hookrightarrow \A^n$,
respectively) we let $\Lin (Y)$ denote the group of all
automorphisms of $Y$ which extend linearly to the ambient space.

\subsection{Automorphisms of smooth projective hypersurfaces}
In the following theorems we gather some results concerning the
groups $\Lin$, $\Aut$, and $\Bir$ for projective hypersurfaces;
see Matsumura and Monsky \cite{MM}, Iskovskikh and Manin
\cite{IM}, Pukhlikov \cite{Pu1}-\cite{Pu3}, Cheltsov \cite{Ch}, de
Fernex, Ein and Must\c{a}t\u{a} \cite{DEM}.

\bthm\label{mamo1} Let $Y$ be a smooth hypersurface in $\PP^n$ of
degree $d$. Then for all $d,n\ge 3$ except for $(d,n)=(4,3)$,
$$\Aut (Y) =\Lin (Y)\,$$ and this group is finite.
It is trivial for a general hypersurface of degree $d\ge 3$.\ethm

There is a similar result for Schubert hypersurfaces in flag
varieties, see Theorem 8.8 in \cite{Te}. For the group of
birational transformations, the following hold.

\bthm\label{mamo2} For $Y\subseteq\PP^n$ as above and for all $d>
n\ge 2$ except for $(d,n)= (3,2)$, $$\Bir (Y)=\Aut(Y)\,.$$ This
group is finite except in the case $(d,n)=(4,3)$ of a smooth
quartic surface $Y\subseteq\PP^3$, where it is discrete, but can
be infinite and different from $\Lin(Y)$ which is finite. The
group $\Bir(Y)$ of a very general quartic surface
$Y\subseteq\PP^3$ is trivial.\ethm

The case $d\le n$ is much more complicated. However, in this case
there are deep partial results, see e.g., \cite{IM, Pu1, DEM}.

Let us indicate briefly some ideas used in the proofs. In case
$d\neq n+1$ the proof is easy and exploits the fact that the
canonical divisor $K_Y=\cO_Y(d-n-1)$ is $\Aut(Y)$-stable. In case
$d=n+1$ the equalities $\Bir(Y)=\Aut (Y)=\Lin(Y)$ follow since
such hypersurfaces represent Mori minimal models. Indeed a
birational map between minimal models is an isomorphism in
codimension 1, see e.g., \cite{KM}, hence it induces an
isomorphism of the corresponding Picard groups. If $n\ge 4$ then
$\Pic(Y)\simeq\Z$ by the Lefschetz Hyperplane Section Theorem.
Therefore any birational transformation $\varphi$ of $Y$ acts
trivially on $\Pic (Y)$ and so preserves the complete linear
system of hyperplane sections $|\cO_Y(1)|$. Since $Y$ is linearly
normal\footnote{I.e., $H^0(\cO_{\PP^n}(1))\to H^0(\cO_Y(1))$ is a
surjection.}, $\varphi$ is induced by a projective linear
transformation of the ambient projective space $\PP^n$. For the
proof of finiteness of the group $\Lin(Y)$ and its triviality for
general hypersurfaces, we refer to the classical paper of
Matsumura and Monsky \cite{MM}.

By virtue of the Noether-Lefschetz Theorem, these arguments can be
equally applied to very general smooth surfaces in $\PP^3$ of
degree $d\ge 4$. For an arbitrary smooth surface $Y$ in $\PP^3$,
the minimality of $Y$ should be combined with the fact that
$\Pic(Y)$ is torsion free. Indeed, any smooth surface in $\PP^3$
of degree $d\ge 4$ represents a minimal model and so is not
birationally ruled, hence its birational automorphisms are
biregular; see e.g., \cite[Theorem 1-8-6]{Mat}. For $d>4$ the
canonical class $K_Y$ defines an equivariant polarization of $Y$,
and $\cO_Y(1)\sim\frac{1}{d-4}K_Y$. Since $\Pic(Y)$ is torsion
free and $H^0(\cO_{\PP^3}(1))\to H^0(\cO_Y(1))$ is a surjection,
all automorphisms of $Y$ are linear. This is not true, in general,
in the case of a smooth quartic surface in $\PP^3$. An example of
such a surface with infinite automorphism group due to Fano and
Severi is discussed in \cite[Theorem 4]{MM}. A non-linear
biregular involution exists on any smooth quartic in $\PP^3$
containing skew lines, for instance, on the Fermat quartic
$x^4+y^4+z^4+u^4=0$; see Takahashi \cite{Ta}.

For a quadric hypersurface $X\subseteq\A^{n+1}$ of dimension $n\ge
2$, the group $\Aut (X)$ is infinite-dimensional \cite[Lemma
1.1]{To}, cf.\ also example \ref{ex} below. For $n=2$ this group
has an amalgamated product structure \cite{DG}; cf.\ also
\cite{ML}.

For a smooth cubic surface $Y\subseteq \PP^3$ the group $\Aut (Y)
=\Lin (Y)$ is finite, while the Cremona group $\Bir (Y)\simeq\Bir
(\PP^2)$ is infinite-dimensional. The automorphism groups of such
surfaces were listed by Hosoh \cite{Ho1} who corrected an earlier
classification by Segre \cite{Se}; see also Manin \cite{Man} and
Dolgachev \cite{Dol1}. The largest order of such a group is $648$.
This upper bound is achieved only for the Fermat cubic surface,
see \cite{Ho2}. The least common multiple of the orders of all
these automorphism groups is $3240=2^3\cdot3^4\cdot 5$ (Gorinov
\cite{Gor}).

The Fermat quartic $x^4+y^4+z^4+u^4=0$ and the smooth quartic
$x^4+y^4+z^4+u^4+12xyzu=0$ in $\PP^3$ can be also characterized in
terms of the orders of their automorphism groups, see Mukai
\cite{Mu}, Kondo \cite{Ko1}, and Oguiso \cite{Og}.

\subsection{Automorphisms of smooth projective varieties}
According to the well known Matsumura Theorem\footnote{Which
generalizes earlier results by Andreotti for surfaces of general
type; see Kobayashi-Ochiai \cite{KO} and Noguchi-Sunada \cite{NS}
for further generalizations.} the group $\Bir (Y)$ of a smooth
projective variety $Y$ of general type is finite, hence also the
group $\Aut(Y)$ is. In particular, this holds if $c_1(Y)< 0$. See
e.g., Xiao \cite{Xi1, Xi2} for effective bounds of orders of
automorphism groups for the general type surfaces. A vast
literature is devoted to automorphism groups of K3 and Enriques
surfaces. These groups are discrete, and infinite in many cases.
Finite groups of automorphism of $K3$-surfaces were classified
e.g., in \cite{Dol2}, \cite{IS}, \cite{Ko2}, \cite{Mas},
\cite{Ni}, \cite{Xi2}.

For varieties of non-general type we have the following result due
to Kalka-Shiffman-Wong \cite{KSW} and Lin \cite[Theorem
II.3.1.2]{LZ}.

\bthm\label{many.auth} Let $Y$ be a smooth projective variety of
dimension $n$. Suppose that not all Chern numbers of $Y$ vanish
and either $c_1(Y)\le 0$ or $H^{n,0} (Y)\neq 0$. Then $\Aut (Y)$
is a discrete group. \ethm

The first assumption is fulfilled, for instance, if $e(Y)=c_n
(Y)\neq 0$, or $e(\cO_Y)\neq 0$, or $c_1^n (Y)\neq 0$, where $e$
stands for the Euler characteristic. However, this assumption does
not hold for an abelian variety $Y=A$. For any projective
embedding $A\hookrightarrow \PP^N$, the group $\Lin (A)$ is
finite, see \cite[\S II.6]{GH}, while the group $\Aut (A)\supseteq
A$ is infinite. Conversely, by Th\'eor\`eme Principale I of
Blanchard \cite{Bl} for any finite subgroup $G\subseteq \Aut (A)$
there exists a projective embedding $A\hookrightarrow \PP^N$ which
linearizes $G$. A general form of Blanchard's Theorem is as
follows (cf.\ \cite[Theorem 3.2.1]{Ak3}).

\bthm\label{blan} Let $Y$ be a smooth projective variety and
$G\subseteq\Aut (Y)$ a subgroup which acts finitely
 on $\Pic (Y)$. Then there is a $G$-equivariant
projective embedding $Y\hookrightarrow\PP^N$. \ethm

Indeed, such an embedding corresponds to a very ample
$G$-invariant divisor class. However, if $G$ acts finitely on
$\Pic (Y)$ then the orbit of any ample class is an ample
$G$-invariant class.

For instance, if $\Pic (Y)$ is discrete and $G$ is connected then
$G$ acts trivially on $\Pic (Y)$. Hence there exists a
$G$-equivariant projective embedding $Y\hookrightarrow\PP^N$.

As another example, consider a smooth Fano variety $Y$ embedded by
a pluri-anticanonical system $\varphi_{|-mK_Y|}:Y
\hookrightarrow\PP^N$ for a suitable $m>0$. The canonical bundle
$K_Y$ being stable under the action of the automorphism group
$\Aut (Y)$ on $Y$, this embedding is equivariant and realizes
$\Aut (Y)$ as a closed subgroup of $\PGL_{N+1}(\C)$. In
particular, this applies to the anticanonical embeddings of del
Pezzo surfaces $Y_d\hookrightarrow\PP^d$ of degree $d\ge 3$. A
rauch description of the automorphisms groups of these surfaces is
as follows, see Proposition 10.1.1 in \cite{Dol1} (cf.\ e.g.,
\cite{De}, \cite{dF}, \cite{Ho3}, \cite{DI1, DI2}, , \cite{BB},
\cite{Bla} for more delicate properties).

\bthm\label{dpaut} Let $Y=Y_d$ be a del Pezzo surface of degree
$d\ge 3$. Then the automorphism group $\Aut (Y)$ acts on the
lattice $Q=(\Z K_Y)^\bot\subseteq\Pic (Y)$ preserving the
intersection form. The image of the corresponding homomorphism
$\rho : \Aut (Y)\to O(Q)$ is contained in the Weyl group $W(Q)$.
The kernel of $\rho$ is trivial for $d\le 5$ and is a connected
linear algebraic group of dimension $2d-10$ for $d\ge 6$. More
precisely, the following hold. \bnum\item For $d\le 5$ the group
$\Aut (Y)$ is finite.
\item For $d\ge 6$ the identity component $\Aut_0\, (Y)=\ker (\rho)$
contains a 2-torus $\T_2\simeq (\C^*)^2$, and $\Aut_0\, (Y)=\T_2$
for $d=6$.
\item For $d\ge 7$ besides the 2-torus $\T_2$ the group $\Aut_0\,
(Y)$ contains a subgroup isomorphic to $\A^2_+=(\C_+)^2$. In
particular, for $d=7$ there are a decomposition
$$\Aut (Y)\simeq (\A^2_+\rtimes\T_2)\rtimes\Z/2\Z\,$$
and a faithfull presentation $\Aut_0 (Y)\hookrightarrow\GL_3(\C)$
with image
$$\left(%
\begin{array}{ccc}
  1 & 0 & * \\
  0 & * & * \\
  0 & 0 & * \\
\end{array}%
\right)=\left(%
\begin{array}{ccc}
  1 & 0 & 0 \\
  0 & * & 0 \\
  0 & 0 & * \\
\end{array}%
\right)\cdot \left(%
\begin{array}{ccc}
  1 & 0 & * \\
  0 & 1 & * \\
  0 & 0 & 1 \\
\end{array}%
\right)\,.$$ \item For $d=8$ either $Y\to \PP^2$ is  a blowup at a
point and then $\Aut (Y)\simeq \GL_2\ltimes\A^2_+\,$, or $Y\simeq
\PP^1\times\PP^1$ and then $\Aut (Y)\simeq
(\PGL_2(\C))^2\rtimes\Z/2\Z$.\item Finally for $d=9$,
$Y\simeq\PP^2$ and $\Aut (Y)\simeq\PGL_3(\C)$. \enum\ethm

\brem\label{vfi} An effective $\A^2_+$-action on a del Pezzo
surface $Y$ of degree $d=7$ can be defined via the locally
nilpotent derivations \be\label{derv} \p_{\alpha,\beta} = \alpha
z\frac{\p}{\p x}+\beta z \frac{\p}{\p y}, \quad\mbox{}\quad
(\alpha,\beta)\in \A^2_+\,.\ee Indeed, the induced $\A^2_+$-action
on $\A^3$: $$(\alpha,\beta).(x,y,z)= (x+\alpha z, y+\beta z,
z)\,$$ descends to an action on $\PP^2$ fixing the line $z=0$
pointwise. The blowup at two points on this line preserves the
action. Likewise one defines an $\A^2_+$-action on $Y$ for $d=8$
or $d=9$. \erem

\subsection{Homogeneous and almost homogeneous
varieties}\label{sect1.3} By the Borel-Remmert Theorem
\cite[3.9]{Ak3} any connected, compact, homogeneous K\"ahler
manifold $V$ is biholomorphic to the product $ \Alb(V)\times Y$ of
the Albanese torus and a (generalized) flag variety $Y=G/P$ (i.e.,
$Y$ is the quotient of a connected semisimple linear algebraic
group by a parabolic subgroup)\footnote{See also \cite{SDS} for a
more general result in the projective case.}. It follows that
every simply connected homogeneous compact K\"ahler manifold is a
flag variety and the same is true for a rational projective
homogeneous  variety (for homogeneous compact complex manifolds
satisfying both conditions this was established by Goto [Got]).
Furthermore, Grauert and Remmert \cite{GR} carried over a result
of Chow \cite{Cho} from abstract algebraic to Moishezon varieties.
Namely, they  proved that a homogeneous Moishezon variety is
projective algebraic. Thus, if such a variety is simply connected
or rational, it is a flag variety.

Every flag variety $G/P$ is a projective rational Fano variety
(see \cite{Sn}). Every ample line bundle $L$ on $G/P$ is very
ample (see e.g., \cite{Ch2}, \cite{Ja}, \cite[\S 3.3.2]{La}, or
\cite[Theorem 7.52]{Te}). The complete linear system $|L|$ defines
a $G$-equivariant embedding $Y\hookrightarrow\PP^n$ with a
projectively normal image \cite[Theorem 1.iii]{RR}.

For a maximal parabolic subgroup $P_{\max}\subseteq G$, the Picard
group $\Pic (G/P_{\max})\cong \Z$ is generated by the class of a
unique Schubert divisorial cycle in $G/P_{\max}$, and this class
is very ample. In the case of a Grassmannian this class gives the
Pl\"ucker embedding. For an arbitrary flag variety $G/P$, its
Picard group $\Pic (G/P)$ is also generated by the classes of the
Schubert divisorial cycles; see e.g., \cite{Ch2} or \cite{Po2}.
The set of maximal parabolic subgroups $P_{\max}$ of $G$ which
contain $P$ is finite. Every Schubert divisor class in $\Pic
(G/P)$ is lifted via a surjection $G/P\to G/P_{\max}$, see e.g.
\cite{LL} or \cite{Sn}. A linear combination of these divisors is
very ample if and only if its coefficients are all positive (see
\cite{Br} for the case of a full flag variety; the general case is
similar \cite[Theorem 7.52]{Te}).

For a description of the automorphism groups of flag varieties see
e.g., \cite[\S 3.3]{Ak3}.

The following ``cone theorem'' describes certain {\it almost
homogeneous} complex varieties. It is due to Akhiezer
\cite[Theorem 3]{Ak1} in algebraic context and to Huckleberry and
E.\ Oeljeklaus \cite{HO1} in analytic one\footnote{The smooth
compact case was done first in \cite{Oe}. See also \cite[Ch.\ 2,
\S 3, Theorem 1]{HO2} for real groups.}.

\bthm\label{cone} Let $X$ be an irreducible reduced complex space
of dimension $\ge 2$. Suppose that a connected complex Lie group
acts by biholomorphic transformations on $X$ with an open orbit
$\Omega\subseteq X$ such that the complement $E=X\setminus\Omega$
is a proper analytic subset with an isolated point, say, $0\in E$.
Then the normalization $\nu\colon\tilde X\to X$ is one-to-one and
$\tilde X$ is biholomorphic to a projective or an affine cone over
a flag variety $G/P$ of some semisimple linear algebraic group $G$
under a certain equivariant projective embedding. The isolated
point $0\in E$ corresponds to the vertex of the cone. In
particular, if $(X,0)$ is smooth then $X\simeq \A^n$ or $X\simeq
\PP^n$.\ethm

Thus the variety $X$ as in the theorem equipped with an
appropriate algebraic structure carries a regular almost
transitive group action. If the initial group is a complex linear
algebraic group, then $G$ is its maximal semisimple subgroup
\cite{Ak1}. Given any $G$-equivariant projective embedding
$\varphi_{|H|}:Y=G/P\hookrightarrow\PP^n$, where $\dim Y\ge 1$,
the affine cone $\affcone_H (Y)$ over the image admits a regular
action  transitive off the vertex of a locally direct product
$\tilde G\cdot \C^*$, with $\C^*$ acting by homotheties, where
$\tilde G\to G$ is a finite group cover.

A similar description exists for the class of quasi-projective
$G$-varieties $X$, where $G$ is a connected linear algebraic group
acting on $X$ with an open orbit $\Omega$, provided that there is
an equivariant completion $\bar X$ of $X$ with disconnected
complement $\bar X\setminus \Omega$ \cite[Theorem 2]{Ak1}. See
also \cite{Ak2} for the case that $\bar X\setminus \Omega$ is a
$G$-orbit  of codimension $1$ in $X$ (in this case it is
connected).

\smallskip

An explicit description of almost homogeneous 2-dimensional affine
cones over smooth projective curves is due to Popov \cite{Po1}
(see also \cite{FZ2} for an alternative proof). We recall that a
Veronese cone $V_d$ is the affine cone over a smooth rational
normal curve $\G_d\subset\PP^d$ i.e., a linearly
non-degenerate\footnote{A projective variety $Y\subseteq\PP^n$ is
linearly non-degenerate if it is not contained in any hyperplane.}
smooth curve in $\PP^d$ of degree $d$. All such curves in $\PP^d$
are  projectively equivalent and rational. For normal
2-dimensional cones, Popov's Theorem can be stated as follows.

\bthm[V.\ Popov]\label{po} Let $X$ be the affine cone over a
smooth projective curve $Y$. If $X$ is normal and admits an
algebraic group action transitive in $X\setminus\{0\}$, then $X$
is a Veronese cone $V_d$ for some $d\ge 1$, and $Y$ is a rational
normal curve $\G_d$. \ethm

Popov \cite{Po1} actually classified all almost homogeneous cones
in dimension $2$ with an isolated singularity (not necessarily
normal). Every such cone possesses a linear $\SL(2,\C)$-action
transitive off the vertex. The group $\Aut (X)$ of a Veronese cone
is infinite dimensional and so cannot be linearized under an
affine embedding; see Section 2.3 below.

\section{Groups acting on affine cones}
\subsection{Linear automorphisms of affine cones} Let us start
with the following result.

\bprop\label{liextnew1} Given two affine cones $X_i=\affcone
(Y_i)\subseteq \A^{n_i+1}$ over smooth, linearly non-degenerate,
projective varieties $Y_i\subseteq \PP^{n_i}$ ($i=1,2$) and an
isomorphism $\varphi:X_1\stackrel{\simeq}{\longrightarrow} X_2$,
the differential $d\varphi(0)$ provides a linear isomorphism
$\Phi:\A^{n_1+1}\stackrel{\simeq}{\longrightarrow} \A^{n_2+1}$
which restricts to an isomorphism
$\Phi|_{X_1}:X_1\stackrel{\simeq}{\longrightarrow} X_2$. In
particular $n_1=n_2$, and $Y_1$ and $Y_2$ are projectively
equivalent. \eprop

\bproof By the linear non-degeneracy assumption
$$T_0X_i=\A^{n_i+1},\quad C_0X_i=X_i,\quad\mbox{and}\quad
\PP(C_0X_i)=Y_i,\quad i=1,2\,,$$ where $T_0X_i$ is the Zariski
tangent space to $X_i$ at the vertex $0\in X_i$,  and $C_0X_i$ is
the tangent cone in $0$ (see e.g., \cite[\S 9.7]{CLS}). Now the
assertion follows since $d\varphi(0)$ provides an isomorphism of
the Zariski tangent spaces and sends the cone $C_0X_1$ onto the
cone $C_0X_2$ \cite[\S 7.3]{Da}. In fact $d\varphi(0)$ lifts to an
isomorphism of blowups
$\Bl_0(X_1)\stackrel{\simeq}{\longrightarrow}\Bl_0(X_2)$
preserving the exceptional divisors. These divisors are isomorphic
to $Y_1$ and $Y_2$, respectively, and $d\varphi(0)$ induces a
linear isomorphism $Y_1\stackrel{\simeq}{\longrightarrow} Y_2$.
\eproof

\brem\label{nrem10} The isomorphism $\varphi$ as in Proposition
\ref{liextnew1} does not need to be linear itself. However, this
is the case under the additional assumption that $Y_1$ is not
birationally ruled (see Proposition \ref{liextnew} below). A
birationally ruled projective variety is a variety birationally
equivalent to a product $Z\times\PP^1$. Recall also that a
birational map $\bar f:\bX_1\dashrightarrow \bX_2$ is said to be
{\em isomorphism in codimension one} if there are subsets
$B_i\subseteq \bX_i$ of codimension at least 2 such that $$\bar
f|(\bX_1\setminus B_1):\bX_1\setminus B_1\to \bX_2\setminus B_2$$
is an isomorphism. \erem

\bprop\label{liextnew} Consider the affine cones
$X_i=\affcone(Y_i)\subseteq \A^{n_i+1}$ over projective varieties
$Y_i\subsetneq \PP^{n_i}$, $i=1,2$. Suppose that $Y_1$ and $Y_2$
are smooth, irreducible, and linearly non-degenerate. If $Y_1$ is
not birationally ruled then every isomorphism
$\varphi:X_1\stackrel{\simeq}{\longrightarrow} X_2$ extends to a
unique linear isomorphism
$\A^{n_1+1}\stackrel{\simeq}{\longrightarrow} \A^{n_2+1}$. In
particular $n_1=n_2$, and $Y_1$ and $Y_2$ are projectively
equivalent. \eprop

The proposition follows immediately from Lemmas \ref{tot} and
\ref{liext3} below. Before passing to the lemmas, let us give two
corollaries, which are the main results of this subsection.

\bcor\label{liext} Let $X=\affcone (Y)\subseteq \A^{n+1}$ be the
affine cone over a smooth projective variety $Y\subseteq \PP^{n}$.
If $Y$ is not birationally ruled then $\Aut\, (X)=\Lin\, (X)$.
Moreover, $\Aut\, (X)$ is a central extension of the group $\Lin\,
(Y)$ by $\C^*$. \ecor

Indeed, the exact sequence
$$0\to \C^*\to \GL(n+1,\C) \to \PGL(n+1,\C)\to 0\,$$ yields
the following one: \be\label{exse234} 0\to \C^*\longrightarrow
\Lin (X)\stackrel{\pi}{\longrightarrow} \Lin (Y)\to 0\,.\ee

\bcor\label{liext2} Let $X=\affcone (Y)$ be the affine cone over a
smooth projective 3-fold $Y$. Suppose that $Y$ is rationally
connected and non-rational. Then $\Aut(X)=\Lin (X)$.\ecor

\bproof Indeed if $Y$ were birationally ruled i.e., birational to
a product $Z\times\PP^1$, then $Z$ would be rationally connected
and so a rational surface. Hence $Y$ would be rational too,
contrary to our assumption. Thus Corollary \ref{liext} applies and
gives the assertion. \eproof

\bexa\label{mfcl} For instance, if $Y\subseteq\PP^n$ is a
non-rational Fano 3-fold and $X=\affcone (Y)$, then $\Aut (X)=\Lin
(X)$. As an example, one can consider any smooth cubic or quartic
3-fold $Y\subseteq \PP^4$.\eexa

For the proof of the next lemma we refer the reader to
\cite[Theorem 2.2]{To} or \cite[Proposition 2.7]{Co}.

\blem\label{tot} Let $X_i=\bX_i\setminus D_i$, $i=1$, $2$, where
$\bX_i$ is a projective variety and $D_i$ an irreducible divisor
on $\bX_i$. Suppose that $X_i$ is regular near $D_i$ for $i=1,2$.
If $D_1$ is not birationally ruled then any isomorphism $f:X_1\to
X_2$ extends to a birational map $\bar f:\bX_1\dashrightarrow
\bX_2$ which is an isomorphism in codimension 1. If in addition
the divisors $D_1$ and $D_2$ are ample then $\bar
f:\bX_1\stackrel{\simeq}{\longrightarrow}\bX_2$ is an isomorphism.
\elem

Proposition \ref{liextnew} is now a direct consequence of the
following lemma.

\blem\label{liext3} Consider two projective varieties
$Y_i\subsetneq \PP^{n_i}$, where $n_i\ge 2$, $i=1$, $2$. Suppose
that $Y_1$ and $Y_2$ are smooth, irreducible, and linearly
non-degenerate. Consider also the affine cones $X_i=\affcone
(Y_i)\subseteq \A^{n_i+1}$  over $Y_i$ and the projective cones
$\bX_i\subseteq \PP^{n_i+1}$, $i=1,2$.
 Let $\varphi:X_1\stackrel{\simeq}{\longrightarrow}
X_2$ be an isomorphism such that the induced birational map
$\bar\varphi:\bX_1\dashrightarrow \bX_2$ is an isomorphism in
codimension 1. Then $\varphi$ extends to a unique linear
isomorphism $\Phi:\A^{n_1+1}\stackrel{\simeq}{\longrightarrow}
\A^{n_2+1}$. In particular $n_1=n_2$, and $Y_1$ and $Y_2$
($\bX_1$ and $\bX_2$, respectively) are projectively
equivalent.\elem

\bproof We let $D_i=\bX_i\setminus X_i$ denote the divisor at
infinity; it is a scheme-theoretic hyperplane section. Since $D_1$
and $D_2$ are ample then (similarly as in Lemma \ref{tot})
$\varphi$ extends to an isomorphism $\bar \varphi:\bX_1\to \bX_2$,
which sends $0\in\A^{n_1+1}$ to $0\in\A^{n_2+1}$. Indeed, these
points are the only singular points of the projective cones
$\bX_1$ and $\bX_2$. Moreover, $\bar \varphi$ sends the generators
of the cone\footnote{That is the projective lines on $\bX_1$
passing through the origin.} $\bX_1$ into generators of $\bX_2$.
Indeed, every generator $l_1$ of $\bX_1$ meets $D_1$ transversally
in one point. The image $l_2=\varphi(l_1)\subseteq \bX_2$
possesses similar properties, hence $l_2$ is again a projective
line through the origin i.e., a generator of the cone $\bX_2$.

It follows that the orbits of the $\C^*$-action on $\bX_1$ are
sent to the orbits of the $\C^*$-action on $\bX_2$. Furthermore
$\bar\varphi$ is $\C^*$-equivariant, hence it induces an
isomorphism $\varphi^*:\cO(Y_2)\stackrel{\simeq}{\longrightarrow}
\cO(Y_1)$ of the homogeneous coordinate rings. These graded rings
are the coordinate rings of the affine cones $X_1$ and $X_2$,
respectively, generated by their first graded
pieces\footnote{Consisting of the restrictions to $X_i$ of linear
functions on $\A^{n_i+1}$, $i=1,2$.}. The graded isomorphism
$\varphi^*$ restricts to a linear isomorphism, say,
$\Psi:\A^{n_2+1}\stackrel{\simeq}{\longrightarrow} \A^{n_1+1}$
between these first graded pieces. The dual isomorphism
$\Phi=\Psi^{\vee}:\A^{n_1+1}\stackrel{\simeq}{\longrightarrow}
\A^{n_2+1}$ provides a desired linear extension of $\varphi$. The
uniqueness of such an extension follows immediately, since $Y_1$
and $Y_2$ are assumed to be linearly non-degenerate.\eproof

For a projective variety $Y\subseteq\PP^n$ with affine cone
$X=\affcone (Y)$ it can happen that $\Aut (Y)\neq \Lin (Y)$, while
$\Aut (X)= \Lin (X)$, as in the following examples.

\bexas\label{abli} 1.\ Let $A$ be an abelian variety. Consider a
projective embedding $A\hookrightarrow\PP^n$ (for instance, a
smooth cubic in $\PP^2$) with affine cone $X=\affcone (A)$. By
Corollary \ref{liext} $\Aut (X)=\Lin (X)$. While $\Lin (A)$ is a
finite group (see \cite[\S II.6]{GH} or Section 1 above),  the
group $\Aut (A)$ contains the subgroup of translations and so is
infinite. Thus  $\Aut (Y)\neq \Lin (Y)$ (cf.\ Blanchard's Theorem
\ref{blan}).

2.\ A smooth quartic $Y\subseteq \PP^3$ is a K$3$-surface and so
is not birationally ruled. Hence again $\Aut (X)=\Lin (X)$, where
$X=\affcone (Y)\subseteq\A^4$. Moreover, $\Lin (Y)$ is a finite
group, while the group $\Aut (Y)$ can be infinite, see the
discussion in \S 1.1. Clearly, non-linear automorphisms of $Y$ are
not induced by automorphisms of $X$.\eexas

\subsection{Lifting $G$-actions to affine cones}
In this subsection we address the following questions.

\smallskip

\noindent \bnum {\em \item When a $G$-action on $Y$ is induced by
a $G$-action on $X$?\item When a $G$-action on $Y$ is induced by a
$\tilde G$-action on $X$? }\enum

A related question is:

\smallskip

{\em Which projective representations can be lifted to linear
ones?}

\smallskip

\noindent Simple examples show that one needs some restrictions on
such a projective representation.
 In the first example below the group  $G$ is finite,
and is connected algebraic in the second.

\bexas\label{yuex} 1.\ The standard representation  on $\A^2$ of
the group of quaternions $Q_8=\{\pm 1, \pm i, \pm j, \pm k\}$
induces a faithful representation of $Q_8$ on any Veronese cone
$V_d\simeq \A^2/\Z_d$ with $d$ odd (cf.\ Subsection \ref{sect2-4}
below). The latter representation descends to an effective linear
action on $\PP^1$ of the dihedral group
$$D_2=Q_8/Z(Q_8)\simeq (\Z/2\Z)^2
\,.$$ However, this $D_2$-action on $\PP^1$ cannot be lifted to a
$D_2$-action on $\A^2$ or on any of the Veronese cones $V_d$ with
$d$ odd. Indeed, otherwise the exact sequence $$0\longrightarrow
Z(Q_8)\longrightarrow Q_8\longrightarrow D_2\longrightarrow 0$$
would split, which is not the case. In other words, the faithful
projective representation $D_2\to\PGL_{n+1}(\C)$ induced by the
Veronese embedding $\varphi_{|\cO_{\PP^1}(n)|}:\PP^1\to\PP^n$
lifts to a linear representation $D_2\to\GL_{n+1}(\C)$ if and only
if $n=2k>0$ is even and so $\cO_{\PP^1}(n)=-k K_{\PP^1}$.

2.\ The standard projective representation of $G=\PGL_2(\C)$ on
$\PP^1$ induces a linear $G$-action on the rational normal curve
$\G_d\subseteq \PP^d$. Suppose that the latter action can be
lifted to the Veronese cone $V_d=\affcone
(\G_d)\subseteq\A^{d+1}$. This would give an irreducible
representation of $G=\PGL_2(\C)$ of dimension $d+1$. However, such
a representation does exist only for $d$ even.  Indeed, every
irreducible representation of $\PGL_2(\C)$ yields an irreducible
representation of $\SL_2(\C)$ trivial on the center, and vice
versa.  \eexas

\brem\label{schur} Concerning question (2), recall that for any
perfect group $G$ there exists a unique universal central
extension (or Schur cover) $G'$ of $G$ such that every projective
representation of $G$ is induced by a linear representation of
$G'$ (see \cite[\S 7]{St}). For a finite perfect group $G$, the
Schur cover $G'$ is again finite. For a perfect (e.g.,
semi-simple) connected linear algebraic group $G$ over $\C$, the
Schur cover  is just the simply connected universal covering group
$G'$ of $G$. \erem

Any linear action $G\to\Lin (X)$ on the affine cone $X=\affcone
(Y)$ induces (via the exact sequence (\ref{exse234})) a linear
action $G\to\Lin (Y)$ on $Y$. The latter factorizes through the
action on $Y$ of the quotient group $G/(G\cap\, \C^*)$. Answering
question (2) above, in the following proposition we provide a
simple criterion as to when a $G$-action on $Y$ is induced by a
$\tilde G$-action on $X=\affcone (Y)$, where $\tilde G$ is a
central extension of $G$ (which is not a Schur cover). The proof
is straightforward. Let us remind that any linear action $G\to
\Lin (Y)$ of a group $G$ on a projective variety $Y\subseteq\PP^n$
stabilizes the very ample divisor class $[\cO_Y(1)]\in\Pic (Y)$.

\bprop \label{imprem} \bnum\item[(a)] Let $Y$ be a smooth
projective variety and $G\to \Aut (Y)$ be a group action on $Y$.
If this action stabilizes a very ample divisor class $|H|\in \Pic
(Y)$, then it extends linearly to the ambient projective space
$\PP^n=\PP H^0 (Y,\cO_Y(H))$.
\item[(b)] Furthermore, let $X=\affcone_H (Y)$ be the affine cone over
$\varphi_{|H|} (Y)$. Consider the central extension $\tilde
G=\pi^{-1}(G)\subseteq \Lin (X)$ of $G$ by $\C^*$, where $\pi
:\Lin (X)\twoheadrightarrow\Lin (Y)$ is as in (\ref{exse234}).
Then the group $\tilde G$ acts linearly on $X$ inducing the given
$G$-action on $Y$.  \enum\eprop

\bcor\label{nnco} Let $G$ be a connected linear algebraic group.
Then any regular $G$-action on a smooth projective variety
$Y\subseteq\PP^n$ is induced by a regular $\tilde G$-action on the
affine cone $\affcone (Y)$, where $\tilde G=\pi^{-1}(G)\subseteq
\Lin (X)$ is a central extension  of $G$ by $\C^*$. \ecor

\bproof Since $G$ is connected, $G$ acts on $\Pic_0 (Y)$. The
group $G$ being a rational variety \cite{Ch1}, every morphism of
$G$ to the abelian variety $\Pic_0 (Y)$ is constant. Hence the
$G$-action on $\Pic_0 (Y)$ is trivial, and so is the induced
action on the Neron-Severi group $NS(Y)=\Pic (Y)/\Pic_0 (Y)$. Thus
$G$ acts trivially on $\Pic (Y)$. By Proposition \ref{imprem}(a)
the $G$-action on $Y$ extends linearly to $\PP^n$. Now the result
follows. \eproof

\brem\label{ratcon} Instead of referring to \cite{Ch1} one can
show directly that every morphism $f:G\to A$ to an abelian variety
$A$ is constant. Clearly, $f$ is constant on any abelian subgroup
of $G$ and on its cosets. Hence $f$ is also constant on any
solvable subgroup. In particular, it is constant on $\Rad(G)$ and
on its cosets. Thus $f$ induces a morphism $G/\Rad(G)\to A$. So we
may assume that $G$ is semisimple. Consider a maximal torus
$\T\subseteq G$ and the collection of its root vectors
$(H_\alpha)_\alpha\subseteq T_eG={\rm lie} (G)$. The subset
$T_e\T\cup (H_\alpha)_\alpha$ consists of the tangent vectors of
algebraic one-parameter subgroups of $G$ and spans the tangent
space $T_eG$. Hence the differential $df(e)$ vanishes. Now the
assertion follows. Indeed, applying left shifts one can produce a
similar situation in any point $g$ of $G$. \erem

A stronger statement holds for pluri-canonical or pluri-anticanonical embeddings.

\bprop\label{akankl} Let $Y$ be a smooth projective variety.
Suppose that for some $m\in\Z$ there is an embedding
$\varphi=\varphi_{|mK_Y|}:Y\hookrightarrow\PP^n$, and let
$X=\affcone (\varphi(Y))$. Then $$\Lin (X) =\C^*\times\Lin
(\varphi(Y))\simeq\C^*\times\Aut (Y)\,,$$ where $\C^*$ acts on the
cone $X$ by scalar matrices.\eprop

\bproof Indeed, the group $\Aut (Y)$ acts on the linear system
$|mK_Y|$ yielding an isomorphism $\Aut (Y)\simeq \Lin
(\varphi(Y))$. Moreover, $\Aut (Y)$ acts on the linear bundle
$\cO(mK_Y)$. Hence it acts linearly on $H^0(Y, \cO(mK_Y))$. The
dual action on $H^0(Y, \cO(mK_Y))^\vee$ preserves  the cone $X$.
This gives an embedding $\Aut (Y)\hookrightarrow\Lin (X)$ and a
splitting of the exact sequence
$$0\to\C^*\to\Lin (X)\to\Lin
(\varphi(Y))\simeq\Aut (Y)\to 0\,.$$ Since the subgroup
$\C^*\subseteq\Lin (X)$ is central, the assertions follow.\eproof

This proposition can be applied to the anticanonical embeddings of
del Pezzo surfaces. In the case where there is a $\C_+$-action on
$X$ the group $\Aut (X)$ is infinite
dimensional. For instance, this is so for the cones over del Pezzo
surfaces of degree $d\ge 4$. For $d\ge 7$ there exists a linear
$\A^2_+$-action on $X$. While for $6\ge d\ge 4$ the group $\Aut (Y)$ is
finite or toric, hence any
$\C_+$-action on $X$ is non-linear;
cf.\ Theorem \ref{main1} in the Introduction and
also Theorems \ref{dpaut} and \ref{cuco}.

\subsection{Groups acting on affine cones}

Similarly as in Proposition \ref{liextnew1}, in the case of a
reductive group action a weaker analog of Corollary \ref{liext}
holds without the assumption of birational non-ruledness.

\blem\label{sile} Suppose that a connected reductive group $G$
acts effectively on the affine cone $X\subseteq \A^{n+1}$ over a
smooth linearly non-degenerate projective variety $Y\subsetneq
\PP^n$. Then there is a faithful representation
$\rho:G\to\GL(n+1,\C)$, which restricts to an effective linear
$G$-action on $X$ inducing a linear action of $G$ on $Y$. \elem

\bproof The vertex $0\in X$ is an isolated singular point of $X$,
hence a fixed point of $G$. Since $G$ is reductive, the induced
representation $\rho$ of $G$ on the Zariski tangent space $T_0X$
is faithful (see e.g., \cite{Ak3} or \cite[Lemma 2.7(b)]{FZ2}) and
descends to $Y$ via the projective representation
$\bar\rho:G\to\PGL(n+1,\C)=\GL(n+1,\C)/\C^*$. \eproof

Let us note that for a non-reductive group action, $\rho$ as above
can be trivial. For instance, this is the case for the
$\C_+$-action $t.(x,y)=(x+ty^2,y)$ on $X=\A^2$.

\smallskip

The following theorem is complementary to Corollary \ref{liext};
 cf.\ also \cite{HO1} for (a).

\bthm\label{1.2} We let $X\subseteq \A^n$ ($n\ge 2$) be the affine
cone over a smooth projective variety $Y\subseteq \PP^{n-1}$.
Suppose that
\begin{enumerate}\item[$\bullet$] The group $\Aut
(Y)$ is finite. \item[$\bullet$] A connected algebraic group $G$
of dimension $\ge 2$ acts effectively on $X$ and contains a
1-dimensional torus $\T\simeq\C^*$ acting on $\A^n$ via scalar
matrices.
\end{enumerate} Then the following hold.
\begin{enumerate}
\item[\rm (a)]
$G$ is a solvable group of rank $1$.
\item[\rm (b)]
There exists an $\A^1$-fibration $\theta:X\to Z$, where $Z$ is an
affine variety equipped with a good $\C^*$-action and $\theta$ is
equivariant with respect to the standard $\C^*$-action on $X$.
Furthermore, $Z$ is normal if $X$ is.
\item[\rm (c)]
$Y$ is uniruled via a family of rational curves parameterized by
$(Z\setminus\{\theta(0)\})/\C^*$.
\end{enumerate}
\ethm

\bproof Consider a Levi decomposition $G=\Rad_u(G)\rtimes L$,
where $L\subseteq G$ is a Levi subgroup (i.e., a maximal connected
reductive subgroup) containing $\T$. By Lemma \ref{sile} the
induced representation $\rho$ of $L$ on the Zariski tangent space
$T_0X$ is faithful. Moreover $\T$ (which acts on $T_0X$ by scalar
matrices) is a central subgroup of $L$. Since the group $\Aut (Y)$
is finite and $L$ is connected, the induced action of the quotient
group $L/\T$ on $Y$ is trivial. Thus $L=\T$ is a maximal torus of
$G=\Rad_u(G)\rtimes \T$, and so $G$ is solvable of rank 1.

By our assumption $\dim_\C (G)\ge 2$. Hence the unipotent radical
$\Rad_u(G)$ is non-trivial and contains a one-parameter subgroup
$U\simeq G_a$. All orbits of $U$ are closed in $X$, and the
one-dimensional orbits are isomorphic to the affine line $\A^1$.
Therefore $X$ is affine uniruled. Its coordinate ring $A=\cO_X$ is
graded by the dual lattice $\T^\vee\simeq \Z$. This grading is
actually positive:
$$A=\bigoplus_{k\ge 0} A_k\,.$$ The infinitesimal generator
$\p$ of the induced $G_a$-action on $A$ is a homogeneous locally
nilpotent derivation of $A$ (see e.g., \cite{Re} or \cite{FZ2}).
The ring of invariants $B=\ker (\p)=A^{G_a}$ is a graded
subalgebra of $A$ with $B_0=A_0=\C$. Therefore the affine variety
$Z=\spec (B)$ is endowed by a $\T$-action with a unique attractive
fixed point $0'=\theta(0)$, where $\theta: X\to Z$ is the orbit
map of the $G_a$-action on $X$. Thus $\theta$ is a
$\T$-equivariant surjection induced by the inclusion $B\subseteq
A$ of graded rings. If $A$ is integrally closed in $\Frac (A)$
then also $B$ is. Indeed, let $Z'=\spec (\bar B)$ be the
normalization of $Z$, where $\bar B$ is the integral closure of
$B$ in $\Frac (A)$. Since $X$ is normal the morphism $X\to Z$
factorizes as $X\to Z'\stackrel{\nu}{\longmapsto} Z$. The locally
nilpotent derivation $\p$ stabilizes  $\bar B$  (see e.g.,
\cite{Sei}, \cite{Vas}, or \cite[Lemma 2.15]{FZ1}) and so the
morphisms $X\to Z'\stackrel{\nu}{\longmapsto} Z$ are equivariant
with respect to the induced $\C_+$-actions.  The $\C_+$-action is
trivial on $Z$, hence also on $Z'$ since $Z'\to Z$ is finite. Thus
$B\subseteq \bar B\subseteq \ker\p= B$, so $B=\bar B$ is normal as
soon as $A$ is.

Since a general one-dimensional orbit of $U\simeq G_a$ in $X$ does
not pass through the vertex $0\in X$ and is not contained in an
orbit closure of $T$ (i.e., in a generator of the cone), there is
a Zariski open subset, say, $\Omega$ of $Y$ covered by the images
of these orbits. Taking Zariski closures yields a family of
rational curves parameterized by $(Z\setminus \{0'\})/\C^*$. Thus
$Y$ is uniruled, as claimed. \eproof

\subsection{Group actions on 2-dimensional affine cones}\label{sect2-4}
The following corollary is immediate from Proposition
\ref{liextnew1}.

\bcor\label{verone} Consider two smooth linearly non-degenerate
curves  $Y_i\subseteq \PP^{n_i}$  ($i=1,2$) of degrees $d_i$, and
let $X_i=\affcone (Y_i)\subseteq\A^{n_i+1}$ be the corresponding
affine cones. Then $X_1\simeq X_2$ if and only if these cones are
linearly isomorphic, if and only if $n_1=n_2$, $d_1=d_2$ and $Y_1$
and $Y_2$ are projectively equivalent.\ecor

Similarly, from Corollary \ref{liext} we deduce the following one.

\bcor\label{verone1} Let $X=\affcone (Y)\subseteq\A^{n+1}$ be the
affine cone over a smooth,  non-rational projective curve
$Y\subseteq\PP^n$. Then $\Aut (X)=\Lin (X)$, and this group is a
central extension of the finite group $\Lin (Y)$ by $\C^*$. \ecor

\brems\label{moduli} 1. However, $\Aut (Y)\neq \Lin (Y)$ for an
elliptic curve  $Y\subseteq\PP^n$, see Example \ref{abli}(1).
Consider further a smooth rational curve $Y\subseteq\PP^n$ of
degree $d>n$. Then $Y$ is neither linearly nor projectively
normal. Indeed, $Y$ is a linear projection of the rational normal
curve $\G_{d}\subseteq\PP^{d}$, and $X=\affcone (Y)$ is a linear
projection of the Veronese cone $V_{d}=\affcone (\G_d)$. The
letter projection gives a normalization of $X$. This is not an
isomorphism as it diminishes the dimension of the Zariski tangent
space at the vertex.

2. The normalizations of the affine cones $X_1$ and $X_2$ over two
smooth rational curves $Y_1$ and $Y_2$, respectively, are
isomorphic if and only if $\deg\, (Y_1)=\deg\, (Y_2)$. While in
general the (non-normal) affine surface $X=\affcone (Y)$ admits
non-trivial equisingular deformations arising from deformations of
the projective embedding $Y\hookrightarrow \PP^n$. For instance,
smooth rational curves $Y$ of type $(1,a)$ on a quadric
$\PP^1\times\PP^1\hookrightarrow\PP^3$ vary in a family of
projective dimension $2a+1$. Hence for any $a\ge 3$ the group
$\PSO(4,\C)$ cannot act transitively on this family.

3. Since any group action on an affine cone $X$ lifts to the
normalization, it is enough to restrict to normal cones. For the
normal Veronese cone $X=V_d\subseteq\A^{d+1}$ over a rational
normal curve $Y=\G_d\subseteq\PP^{d}$ we have $\Aut (X)\neq\Lin
(X)$. Moreover, $V_d\simeq\A^2/\Z_d$ being a toric surface, for
every $d\ge 1$ the group $\Aut (V_d)$ is infinite dimensional. In
particular this is not an algebraic group. Indeed, the graded
coordinate ring $\cO(V_d)$ admits a nonzero locally nilpotent
derivation $\p$ corresponding to an effective $\C_+$-action on
$V_d$ \cite{FZ2}. The kernel $\ker (\p)\subseteq \cO(V_d)$ is
isomorphic to the polynomial ring $\C[t]$. For any $p\in\C[t]$,
the derivation $p\cdot\p$ is again locally nilpotent. Thus
$\C[t]\cdot\p$ is the Lie algebra of an infinite dimensional
abelian subgroup $G\subseteq \Aut (V_d)$.

4. There are actually two independent $\C_+$-actions on $V_d$ with
different orbits, and even a continuous family of such actions;
see e.g., \cite{FZ2}. Danilov and Gizatullin \cite{DG} studied the
structure of an amalgamated product on the group $\Aut (V_d)$,
while Makar-Limanov  \cite{ML} provided an explicit description of
this group.

5. Similarly, independent $\C_+$-actions, and an amalgamated
product structure, exist on any normal affine toric surface
different from $\A^1_*\times\A^1$ or $\A^1_*\times\A^1_*$. Every
such surface is of the form $V_{d,e}=\A^2/\Z_d$ for an appropriate
diagonal action
$$\zeta . (x,y)=(\zeta x, \zeta^e y), \quad\mbox{where}\quad
\zeta^d=1\quad\mbox{and}\quad \gcd(e,d)=1,$$ of the cyclic group
$\Z_d=\langle\zeta\rangle$ on the affine plane $\A^2$. Choosing a
system of homogeneous generators in the graded coordinate ring
$\cO(V_{d,e})$ yields an embedding $V_{d,e}\hookrightarrow
\A^{N+1}$, which is equivariant with respect to a suitable
diagonal $\C^*$-action on $\A^{n+1}$ with positive weights
$w=(w_0,\ldots,w_N)$. In this way $V_{d,e}$ can be realized as the
affine cone over a smooth rational curve in the corresponding
weighted projective space $\PP^N_w=(\A^{N+1}\setminus\{0\})/\C^*$.

6. By Popov's Theorem \cite{Po1}\footnote{This is actually an
earlier version of the Cone Theorem in dimension 2; see Section
1.} the Veronese cones $V_d=V_{d,1}$ can be characterized as
normal affine surfaces on which an algebraic group acts with an
open orbit and a fixed point. \erems

Let us construct an explicit example of a non-linear biregular
automorphism of a Veronese cone $V_d$   for every $d\ge 1$.

\bexa\label{vero} Consider as before the Veronese cone
$V_d\subseteq\A^{d+1}$ over a rational normal curve
$\G_d\subseteq\PP^{d}$. Let $\bar V_d\subset\PP^{d+1}$ be the
Zariski closure, and $\tilde V_d$ be the blowup of $\bar V_d$ at
the vertex $0\in V_d$. It is well known \cite{DG} that $\tilde
V_d\simeq\Sigma_d$, where $\Sigma_d$ denotes a Hirzebruch surface
with the exceptional section $S_0$ and a disjoint section at
infinity, say, $S_\infty$ with $S_0^2=-d$ and $S_\infty^2=d$.
Therefore
$$V_d\simeq (\Sigma_d\setminus S_\infty)/S_0\,.$$ To exhibit a
non-linear automorphism of $ V_d$ is the same as to exhibit an
automorphism of $\Sigma_d\setminus S_\infty$, which extends to a
birational transformation of $\Sigma_d$ preserving the exceptional
section $S_0$ but not the ruling $\pi:\Sigma_d\to\PP^1$ (or,
equivalently, which blows down the curve $S_\infty$). On the level
of dual graphs, such a birational transformation consists e.g., in
the following sequence of blowups and blowdowns \cite{FKZ}:

$$\bdi
\raisebox{1mm} {$ \cou{d}{S_\infty}$ } \quad &\rDashto & \quad
\raisebox{1mm} {$\cou{-1}{S_\infty}
\lin\!\!\!\!\!\!\lin\cou{-1}{v_{d+1}}\lin \boxo{A_d}$}
 \raisebox{1mm} \quad & \rDashto & \quad
\raisebox{1mm} {$\cou{0}{v_{d+1}}\lin \boxo{A_d}$} \quad &
\rDashto & \quad \raisebox{1mm} {$\cou{-1}{S_\infty'}\lin
\cou{-1}{v_{d+1}}\lin \boxo{A_d}$} \quad &\rDashto & \quad
\raisebox{1mm} {$\cou{d}{S'_\infty}$}\quad. \edi
$$

\smallskip

\noindent Here a box marked $A_d$ represents the linear chain
$[[-2,\ldots,-2]]$ of length $d$. The centers of blowups on the
curves $S_\infty$ and $v_{d+1}$ can vary. Anyhow, the section
$S_\infty$ being contracted, the resulting biregular
transformation of the Veronese cone $V_d$ is non-linear. \eexa

\section{Group actions on 3-dimensional affine cones} The main result of this
section is the existence of a $\C_+$-actions on the affine cones
over every smooth del Pezzo surface of degree $\ge 4$. The proof
exploits a general geometric criterion for the existence of such
an action.

\subsection{Existence of $\C_+$-actions on affine cones: a geometric criterion}
\bsit\label{genco} Let $Y$ be a smooth projective variety, and let
$H\in\Div (Y)$ be an ample polarization of $Y$. Consider the total
space $\hat X$ of the line bundle $\cO_Y(H)$ with the zero section
$S_0\subseteq \hat X$. Under the natural identification $S_0\simeq
Y$, we have $\cO_{S_0}(S_0)=\cO_{Y}(-H)$. Hence $S_0$ is
contractible, i.e., there is a birational contraction $\upsilon
:\hat X\to X$, where $X$ is a normal affine variety and
$\upsilon(S_0)$ is a point. In this situation, we call $X$ a
\textit{generalized cone} over $(Y,H)$. If $H$ is very ample, then
$X$ coincides with the normalization of the usual affine $\affcone
(Y)$ cone over $Y\hookrightarrow \PP^{n}$, where the embedding is
given by the linear system $|H|$. So we write
$X=\affcone_H(Y)_{\rm norm}$. In this section we provide a
criterion of existence of a $\C_+$-action on a generalized cone.

Let us note that $X$ can be compactified to the projective cone
$\bar X$ over $Y$ by adding a divisor at infinity $S_\infty\simeq
Y$. The divisor $S_\infty$ on $\bar X$ being ample, the variety
$X=\bar X\setminus S_\infty$ is affine. \esit

\bsit For instance, the affine cone over $\PP^2$ in $\A^3$
coincides with $\A^3$ and so admits a transitive action of the
additive group $\C_+^3$. In the following example we exhibit an
effective $\C_+^2$-action on the affine cone $X\subseteq\A^4$ over
a smooth quadric $Y\subseteq\PP^3$ (cf.\ another constructions in
\cite{Sh}). The automorphism groups of affine quadrics were
studied e.g., in \cite{DG, Doe, To}. Over a general base field,
this group is infinite dimensional as soon as the corresponding
quadratic form is isotropic \cite[Lemma 1.1]{To}. The proof of
Lemma 1.1 in \cite{To} provides a nontrivial linear $\C_+$-action
on any quadric over $\C$.  In the following example we exhibit an
explicit effective $\C_+^2$-action on the affine cone over a
smooth quadric in $\PP^3$.\esit

\bexa\label{ex} All smooth quadrics in $\PP^3$ are projectively
equivalent. Choosing for instance the quadric
$$Y=\{xy=zu\}\,$$ we can define
a linear $\C_+^2$-action on $X={\rm AffCone}\, (Y)\subseteq\A^4$
by the following pair of commuting locally nilpotent derivations
on the ring $A=\cO(X)$:
$$\partial_1=u\p/\p x+y\p/\p z\quad\mbox{and}\quad\partial_2=
u\p/\p y+x\p/\p z\,.$$ See \cite{AS, Sh} for a more thorough
treatment on the subject. \eexa

Let us introduce the following notion.

\bdefi\label{cyl} Let $X$ be an affine variety. For a function
$f\in \cO(X)$ we let $$\DD_+(f)=X\setminus
\V_+(f),\qquad\mbox{where}\quad \V_+(f):=f^{-1}(0)\,.$$ We say
that $X$ is {\em cylindrical} if $X$ contains a dense principal
Zariski open subset $U=\DD_+(f)$ isomorphic to the cylinder
$Z\times\A^1$ over an affine variety $Z$. \edefi

The following proposition generalizes Lemma 1.6 in \cite{FZ2}.

\bprop\label{ext} For an irreducible affine variety $X$, the
following conditions are equivalent: \bnum\item[(i)] $X$ possesses
an effective $\C_+$-action.\item[(ii)] $X$ is cylindrical. \enum
\eprop

\bproof First we suppose that $X$ possesses an effective
$\C_+$-action $\psi$ with the associate locally nilpotent
derivation $\p\neq 0$. The filtration \be\label{filt}
0\in\ker\p\varsubsetneq \ker \p^{(2)}\varsubsetneq\ker
\p^{(3)}\ldots\ee being strictly increasing, we can find
$g\in\cO(X)$ such that $\p^{(2)} g=0$ but $h:=\p g\neq 0$. Thus
$\p h=0$ and so $h\in\cO(X)$ is $\psi$-invariant. Letting $s=g/h$
and $U=\DD_+(h)$ the function $s\in \cO(U)$ gives a {\em slice} of
$\p$ that is, $\p (s)=1$. Consequently, the restriction of $s$ to
any 1-dimensional orbit $O$ of $\psi$ in $U$ is an affine
coordinate on $O\simeq\A^1$. By the Slice Theorem (\cite[Cor.\
1.22]{Fr}), $\cO(U)\simeq \ker (\p)[s]$ and $\p=\p/\p s$.
Therefore $U\simeq Z\times\A^1$, where $Z=\Spec (\ker \p)\simeq
s^{-1}(0)$. This yields (ii).

To show the converse, assume that $X$ is cylindrical. Let
$U=\DD_+(f)\simeq Z\times\A^1$ be a principal cylinder in $X$ as
in Definition \ref{cyl}. We consider the natural $\C_+$-action
$\phi$ on $U$ by translations along the second factor. Since $f|U$
does not vanish it is constant along any orbit of $\phi$ and so
$\phi$-invariant. Letting $\p$ denote the locally nilpotent
derivation on $\cO(U)$ associated to $\phi$, the derivation
$\p_n:=f^n\p\in\Der (\cO(U))$ is again locally nilpotent for any
$n\in\N$. Let $a_1,\ldots,a_k$ be a system of generators of
$\cO(X)$, and let $N\in\N$ be sufficiently large so that $f^N\p
a_i\in \cO(X)$ for any $i=1,\ldots,k$. Then $\p_N (a_i)\in\cO(X)$
for any $i=1,\ldots,k$, hence $\p_N (\cO(X))\subseteq\cO(X)$. Thus
the derivation $\p_N|{\cO(X)}\in\Der (\cO(X))$ is locally
nilpotent and so generates an effective $\C_+$-action $\psi$ on
$X$. Therefore (i) holds. \eproof

\brem\label{smo} Clearly the $\C_+$-actions $\phi$ and $\psi|U$ as
in the proof have the same orbits, and $\V_+(f)=\{f=0\}$ consists
of fixed points of $\psi$. \erem

In the case of  affine cones, Theorem \ref{ext1} below gives a
more practical criterion. We need the following definition.

\bdefi\label{polar} For a projective variety $Y$ with a (very
ample) polarization $\varphi_{|H|}:Y\hookrightarrow \PP^n$, we
call an {\em $H$-polar subset} any Zariski open subset of the form
$U=Y\setminus \supp D$, where $D\in |dH|$ is an effective divisor
on $\PP^n$.\edefi

\bsit\label{kami} Recall that an affine ruling on a variety $U$ is
a morphism $\pi:U\to Z$ such that every scheme theoretic fiber of
$\pi$ is isomorphic to the affine line $\A^1$. By a theorem of
Kambayashi and Miyanishi \cite{KaMi} (see also \cite{KaWr,RS,Du}),
every affine ruling $\pi:U\to Z$ on a normal variety $U$ over a
normal base $Z$ is a locally trivial $\A^1$-bundle.\esit

\bthm\label{ext1} Let $Y$ be a smooth projective variety with a
very ample polarization $\varphi_{|H|}:Y\hookrightarrow \PP^n$.
Then the following hold.
\begin{enumerate}\item[\rm (a)]
If the affine cone $X=\affcone_H (Y)$ admits an effective
$\C_+$-action, then $Y$ possesses an $H$-polar open subset $U$,
which is the total space of a line bundle $U\to Z$.
\item[\rm (b)] Conversely, if $Y$ possesses an $H$-polar
open subset $U$ equipped with an affine ruling
$U\stackrel{\A^1}{\longrightarrow} Z$ and a section $Z\to U$,
where $Z$ is smooth and $\Pic (Z)=0$, then the affine cone
$X=\affcone_H (Y)$ admits an effective $\C_+$-action.
\end{enumerate}\ethm

\begin{proof}
(a) Let $\psi'$ be an effective $\C_+$-action on $X$ with
associate locally nilpotent derivation $\p'\neq 0$. Using the
natural grading of the coordinate ring
$$A=\cO(X)=\bigoplus_{i\ge 0} A_i,$$ $\p'$ can be decomposed into a
finite sum of homogeneous derivations $\p'=\sum_{i=1}^n \p'_i$,
where the principal component $\p:=\p'_n\neq 0$ is again locally
nilpotent. The $\C_+$-action $\psi$ on $X$ generated by $\p$
extends to an effective action of a semi-direct product $G=\C_+
\rtimes\,\C^*$ on $X$.

The filtration \eqref{filt} from the proof of Proposition
\ref{ext} consists now of graded subrings. Hence we can find
homogeneous elements $\hat g,\,\hat h\in A$ such that $\p \hat
g=\hat h$ and $\p\hat h=0$. In the notation of \ref{cyl} we let
$$\hat U=\DD_+(\hat h)\subseteq X\quad\mbox{and}\quad \hat
Z=\V_+(\hat g)\setminus\V_+(\hat h)\subseteq\hat U\,.$$ Likewise
 in the proof of Proposition \ref{ext}, we obtain a
decomposition $\hat U\simeq \hat Z \times\A^1$.

Furthermore, $G$ acts on $\hat U\simeq \hat Z\times\A^1$
respecting the product structure. More precisely, $\C_+$ acts by
shifts on the second factor i.e., along the fibers of the morphism
$\hat\pi:\hat U\to\hat Z$. Since $\hat g$, $\hat h$, and $\p$ are
homogeneous, $\C^*$ acts on $\hat U$ stabilizing $\hat Z$ and
sending the fibers of $\hat\pi$ into fibers. The factorization by
the $\C^*$-action on $\hat U$ yields a Zariski open subset $U=\hat
U/\C^*\subseteq Y$ and a divisor $Z=\hat Z/\C^*$ on $U$ so that
\be\label{minuszero} \hat U=\affcone (U)\setminus
\{0\}\quad\mbox{and}\quad \hat Z=\affcone (Z)\setminus \{0\}\,.\ee
The map $\hat\pi$ defines an affine ruling
$\pi:U\stackrel{\A^1}{\longrightarrow} Z$ with a section
$Z\hookrightarrow U$. Each fiber of $\pi$ is the quotient of a
$G$-orbit in $U$ by the $\C^*$-action. Since $Y$ and $Z$ are
smooth, by the Kambayashi-Miyanishi Theorem cited in \ref{kami},
$\pi$ is an $\A^1$-bundle and, moreover, a vector bundle since it
possesses a section.

Finally, since $D:=\hat h^*(0)\in |dH|$, where $d=\deg (\hat h)$,
the open set $U\subseteq Y$ is $H$-polar (see Definition
\ref{polar}). This shows (a).

To show (b), suppose that $Y$ possesses an $H$-polar open subset
$U$ with an affine ruling $\pi:U\stackrel{\A^1}{\longrightarrow}
Z$ and a section $Z\to U$, where $Z$ is smooth and $\Pic (Z)=0$.
Since both $U$ and $Z$ are smooth, $\pi$ is locally trivial  by
the Kambayashi-Miyanishi Theorem. Since $\pi$ has a section $Z$,
$\pi:U\to Z$ is a line bundle and $Z$ is the zero section. This
bundle is trivial since $\Pic (Z)=0$. Thus $U\simeq Z\times\A^1$.
In particular $\Pic (U)=0$.

Let further $\sigma: \tilde X\to X$ be the blowup of the vertex
$0\in X$. The induced morphism $\rho :\tilde X\to Y$ has a natural
structure of a line bundle with the exceptional divisor
$E=\sigma^{-1}(0)$ as the zero section. Since $\Pic (U)=0$, the
restriction $\rho|\tilde U: \tilde U\to U$ to $\tilde
U:=\rho^{-1}(U)\subseteq\tilde X$ yields a trivial line bundle.
Hence
$$\tilde U\simeq U\times\A^1\simeq Z\times\A^1\times\A^1\quad\mbox{and,
similarly,}\quad\tilde Z:=\rho^{-1}(Z)\simeq Z\times\A^1\,.$$
Under this isomorphism $E\cap \tilde U$ is sent to
$U\times\{0\}\simeq Z\times\A^1\times\{0\}$ and $E\cap \tilde Z$
to $Z\times \{0\}$. For $\hat U=\tilde U\setminus E$ and $\hat
Z=\tilde Z\setminus E$ as in (\ref{minuszero}) we obtain
\be\label{decpr} \hat U\simeq
Z\times\A^1\times\C^*\quad\mbox{and}\quad \hat Z\simeq Z\times
\{0\}\times\C^*\,.\ee Thus $\hat U\simeq \hat Z\times \A^1$ is a
cylinder in $X$. Since $U\subseteq Y$ is $H$-polar, $\hat
U\subseteq X$ is a principal Zariski open subset and so $X$ is
cylindrical. Now (b) follows by Proposition \ref{ext}.
\end{proof}

\brems\label{genc} 1. This theorem, with the same proof, holds
also for generalized cones (see \ref{genco}). In particular, we
may assume that $H$ is just an ample divisor.

2. It is easily seen that if a cone $X=\affcone_H (Y)$ admits an
effective $\C_+$-action, then also the cone $X_k=\affcone_{kH}
(Y)$ admits such an action for any $k\ge 1$. Moreover, this cone
$X_k$ is normal for $k\gg 1$, see \cite[Ch.\ II, Ex\. 5.14]{Ha}.
\erems

\brem\label{genc1} The construction of a $\C_+$-action on $X$ as
in the proof of (b) can be made more explicite. The product
$G=\C_+\times\C^*$ acts on $\hat U$ preserving the product
structure in (\ref{decpr}):
$$G\ni (a,\lambda):\,\hat U\to\hat U,\quad (z,x,y)
\longmapsto (z,x+a,\lambda y)\,.$$ The generators $\p/\p x$ and
$y\p/\p y$ of the $\C_+$- and $\C^*$-actions commute. Letting
$D=X\setminus \hat U$, there is a regular function
$f\in\mathcal{O} (X)$ such that $\divis (f)=nD$. Moreover, we can
choose $f$ of the form $f=y^kg(z)$, where $g\neq 0$. For $N\gg 1$
the $\C_+$-action generated by $\p = f^N\p/\p x$ extends to the
cone $X$, see the proof of Proposition \ref{ext}. With this new
$\C_+$-action, a semidirect product $\C_+ \rtimes\,\C^*$ acts
effectively on $X$. However, the factors do not commute any more.
\erem

Theorem \ref{ext1} yields the following criterion of existence of
a $\C_+$-action on certain 3-dimensional affine cones.

\bcor\label{ext2} Let $Y$ be a  rational smooth projective surface
with a polarization $\varphi_{|H|}:Y\hookrightarrow \PP^n$, and
let $X=\affcone_H (Y)\subseteq\A^{n+1}$ be the affine cone over
$Y$. Then $X$ admits a nontrivial $\C_+$-action if and only if $Y$
possesses an $H$-polar cylinder $U\simeq Z\times \A^1$, where $Z$
is a smooth affine curve.\ecor

\begin{proof} Since $Y$ is smooth and
rational, $Z$ as in Theorem \ref{ext1} is a  non-complete smooth
rational curve. Thus $\Pic (Z)=0$. Hence the affine rulings from
Theorem \ref{ext1} and its proof are actually direct products. Our
assertion can be easily deduced now from Theorem \ref{ext1}.
\end{proof}

Using this criterion, we show next that for an arbitrary smooth
rational surface $Y$, some affine cone over $Y$ admits a
nontrivial $\C_+$-action.

\bprop\label{ratsurfcon} Let $Y$ be a rational smooth projective
surface. Then there is an embedding $\varphi:Y\hookrightarrow
\PP^n$ such that the affine cone $X=\affcone(\varphi(Y))$ is
normal and admits an effective $\C_+$-action. \eprop

\bproof Any point $Q\in Y$ possesses an affine neighborhood
$U\simeq \A^2$. An argument from \cite[(2.5)]{Fu} shows that
$Y\setminus U$ supports an ample divisor. Indeed, $\Pic (Y)$ is a
free abelian group generated by the components $\Delta_i$ of the
divisor $Y\setminus U$. Hence $\Delta_j^2>0$ for some $j$.

Choose a nef and big effective divisor $D=\sum \delta_i\Delta_i$
such that $D\cdot \Delta_i>0$ whenever $\delta_i>0$, with a
maximal possible value of $\lambda(D):={\rm card} \{i \mid
\delta_i> 0\}$. Assume on the contrary that $\supp (D) \neq \supp
(\sum \Delta_i)$, i.e., $\delta_i=0$ for some $i$. Since $\supp
(\sum \Delta_i)$ is connected, there is a component
$\Delta_k\not\subseteq \supp (D)$ with $D\cdot \Delta_k>0$. Then
for $t\gg 0$ the divisor $tD+\Delta_k$ is again nef and big. This
contradicts our maximality assumption for $\lambda(D)$. Therefore
$\supp (D) = \supp (\sum \Delta_i)$ is ample. So for $m\gg 1$ the
linear system $|mD|$ gives an embedding $Y\hookrightarrow \PP^n$
with a projectively normal image, see Exercise 5.14 in \cite[Ch.\
II]{Ha}. Since $Y$ admits an $|mD|$-polar cylinder, $X$ is normal
and cylindrical. By Corollary \ref{ext2}, $X$ admits an effective
$\C_+$-action, as required.\eproof

The following question arises.

\bsit\label{qstn} {\bf Question.} {\em Does there exist a
polarized smooth rational surface $(Y,H)$ without any $H$-polar
cylinder?} \esit

\brem\label{pryu} If $U$ is an $H$-polar cylinder on $Y$ then it
is also $kH$-polar for any $k\in \N$, and vice versa. Thus the
existence of an $H$-polar cylinder depends only on the ray of $H$
in the ample cone of $Y$. Moreover, since the irreducible
components of the divisor $D=Y\setminus U$ span the Picard group
$\Pic (Y)$ and the ample cone is open, the property of a cylinder
$U$ to be $H$-polar is stable under small perturbation of $H$.
\erem

For any  smooth rational projective surface, the $H$-polar
cylinder from Proposition \ref{ratsurfcon} can be chosen to be
isomorphic to the affine plane. Let us provide similar examples in
higher dimensions.

\bexa\label{flvars} Consider a flag variety $G/P$ with an ample
polarization $H$ (see \S \ref{sect1.3}). By Corollary \ref{nnco}
the $G$-action on $G/P$ lifts to a $\tilde G$-action on the cone
$\affcone_H(G/P)$, where $\tilde G$ is the universal cover of $G$.
The actions of one-parameter unipotent subgroups of $\tilde G$
yield effective $\C_+$-actions on the cone. Actually $G/P$
contains an $H$-polar open cylinder $U$ isomorphic to an affine
space $\A^n$ (cf.\ Theorem \ref{ext1}(a)).

Indeed, let $B_+\subseteq P$ be a Borel subgroup of $G$, and let
$B_-$ be the opposite Borel subgroup so that $B_+\cap B_-$ is a
Cartan subgroup. Then $B_-\cdot P$ is open in $G$ and so the
$B_-$-orbit $U$ of $e\cdot P$ is open in $G/P$. Thus $U$ is a big
Schubert cell. Since $U$ is also an orbit of the maximal unipotent
subgroup $B_u\subseteq B_-$, it is isomorphic to $\A^n$. In
particular, $U$ is a cylinder in $Y$. Letting $D=Y\setminus
U=\bigcup_i D_i$, the Schubert divisors $D_i$ form a basis in
$\Pic (G/P)$. In this basis $H=\sum_i \alpha_iD_i$, where
$\alpha_i>0$ for all $i$ since the divisor $H$ is ample, see
\cite{Sn} or \cite[Theorem 7.53]{Te}. Hence $U$ is an $H$-polar
cylinder in $G/P$.

We note that the action of $\tilde G$ on the affine cone $\tilde
X:=\affcone_H(G/P)$ is transitive off the vertex $0\in \tilde X$.
Indeed, we may suppose that $X=\tilde G/\tilde P$, where $\tilde
G$ is semisimple, simply connected, and $\tilde P\subseteq \tilde
G$ is parabolic. Since $\tilde X$ is affine, the stabiliser
Stab$_{\tilde G} (x)$ of a point $x\in\tilde X\setminus\{0\}$
cannot contain a parabolic subgroup.
Hence the stabilizer Stab$_{\tilde G} ([x])$ (conjugate to $\tilde P$)
acts
non-trivially on the generator of the cone through $x$. By
Corollary 1.5 in \cite{Po4} the Makar-Limanov invariant $ML(\tilde
X)$ is trivial (cf. Theorem \ref{mli};
see Section 3.3 below for the definition of the
Makar-Limanov invariant). \eexa

\brem\label{200}  The existence of a cylinder in a projective
variety isomorphic to an affine space is rather exceptional. For
instance, none of smooth rational cubic $4$-folds in $\PP^5$, and
none of smooth $3$-fold intersections of a pair of quadrics in
$\PP^5$ contains a Zariski open set isomorphic to an affine space,
see \cite{PS} and \cite{Pr2}. At the same time, every smooth
intersection $Y$ of a pair of quadrics in $\PP^5$ contains a
$(-K_Y)$-polar cylinder, see Proposition \ref{quaint} below and
its proof. \erem

\subsection{$\C_+$-actions on affine cones over del Pezzo
surfaces} Let us explain the reason why are we interested in the
affine cones over del Pezzo surfaces.

\bsit\label{expli} A normal variety $X$ is
\textit{$\Q$-Gorenstein} if some multiple $nK_X$ of the canonical
Weil divisor $K_X$ is Cartier. This notion is important in the
Mori minimal model program (MMP). It is easily seen that the
generalized cone $X=\affcone_H (Y)$ over a smooth polarized
variety $(Y,H)$ is $\Q$-Gorenstein if and only if $aH\thicksim
-bK_Y$ for some $a\in \N$, $b\in \Z$ (\cite[Example 3.8]{Kol1}).
(If, moreover, $H\thicksim -K_Y$, then $X$ is Gorenstein and has
at most canonical singularity at the origin.)

On the other hand, if $\Aut(X)\neq \Lin(X)$ then by Corollary
\ref{liext} $Y$ is birationally ruled, hence the Kodaira dimension
of $Y$ is negative, see \cite{Kol2}. Thus $b>0$, i.e, $-K_Y$ is
ample. Consequently, $Y$ is a Fano variety.

Therefore, if the affine cone $X$ over $(Y,H)$ is $\Q$-Gorenstein
and admits an effective non-linear $\C_+$-action, then  $Y$ is a
Fano variety and $H\in\Q_{>0}[-K_Y]$. In particular, if $\dim
(Y)=2$ then $Y$ is a del Pezzo surface with its
pluri-anticanonical embedding.\footnote{Cf.\ Remark
\ref{abli}(3).}

From now on we assume that $Y$ is a del Pezzo surface of degree
$d\ge 3$ and $H=-K_Y$ is the anti-canonical polarization.
 Thus the linear system
$|-K_Y|$ is very ample and provides an embedding $Y\hookrightarrow
\PP^{d}$ onto a projectively normal smooth surface of degree $d$,
see e.g., \cite{Dol1}. The affine cone $X=\affcone_{-K_Y} (Y)$ has
a normal, canonical, Gorenstein (hence also Cohen-Macaulay)
singularity at the vertex.\esit

The following theorem is the main result of this subsection (see
Theorem \ref{main1} in the Introduction).

\bthm\label{cuco} Let $Y_d$ be a smooth del Pezzo surface of
degree $d$ anticanonically embedded into $\PP^d$, where $4\le d\le
9$, and let $X_d\subseteq\A^{d+1}$ be the affine cone over $Y_d$.
Then $X_d$ admits a nontrivial $\C_+$-action. \ethm

\bproof Consider a pencil $\cL_{\PP^2}=\langle C_1, C_2\rangle$ on
$\PP^2$ generated by a smooth conic $C_1$ and a double line
$C_2=2l$, where $l$ is tangent to $C_1$ at a point $P_0\in C_1$.
Then $L\setminus\{P_0\}\simeq\A^1$, where $L$ is a general member
of $\cL_{\PP^2}$. Moreover, $U=\PP^2\setminus (C_1\cup C_2)$ is a
cylinder over $\A^1_*$. Blowing up at $9-d$ distinct points $Q_i$
on $C_1\setminus\{P_0\}$, where $9\ge d\ge 4$, we obtain a del
Pezzo surface $Y$ of degree $d$ with a contraction
$\sigma:Y\to\PP^2$, and any such surface can be obtained in this
way, except for $\PP^1\times\PP^1$.

The cylinder $U'=\sigma^{-1}(U)\simeq U$ is $(-K_Y)$-polar (see
Definition \ref{polar}). Indeed, let $E_i=\sigma^{-1}(Q_i)$,
$i=1,\ldots,9-d$. For any $1\gg \varepsilon
>0$ we have
$$-K_{\PP^2}\equiv (1+\varepsilon) C_1 + (1-2\varepsilon)l\,.$$
Hence,
$$-K_{Y}=\sigma^*(-K_{\PP^2})-\sum_{i=1}^{9-d} E_i
\equiv (1+\varepsilon) \Delta_1 +
(1-2\varepsilon)\Delta_2+\varepsilon\sum_{i=1}^{9-d} E_i\,,$$
where $\Delta_1$ and $\Delta_2$ are the proper transforms in $Y$
of $C_1$ and $l$, respectively. Thus $U'$ is a $(-K_Y)$-polar
cylinder on $Y$.

In the remaining case where $Y=\PP^1\times\PP^1$, the natural
embedding $\A^2=\A^1\times\A^1$ into $Y$ yields a $(-K_Y)$-polar
cylinder on $Y$. Applying now Corollary \ref{ext2} ends the proof.
\eproof

The proof exploits a $(-K_{\PP^2})$-polar cylinder on $\PP^2$ made
of a pencil of conics with a common tangent line. Based on the
same idea, we give below some alternative constructions of polar
cylinders on anticanonically polarized del Pezzo surfaces of
degrees $\ge 4$. Due to Corollary \ref{ext2}, this leads to new
$\C_+$-actions on the cones over del Pezzo surfaces of degree $\ge
4$ under their anticanonical embeddings. These examples will be
useful in the sequel.

\bexa\label{misuexa} Consider a pencil of rational curves on
$\PP^2$ with a unique base point $P$. (Similarly, one can find
such a pencil on the quadric $\PP^1\times\PP^1$.) Then the
complement of the union of its degenerate members (or of a general
one, if all members are non-degenerate) is a $(-K_{\PP^2})$-polar
cylinder on $\PP^2$. In \cite{MiSu} an example  was proposed of
such a pencil of quintic curves. Moreover, there is a smooth conic
$C_1$ and a rational unicuspidal quintic $C_2$ from the pencil as
in \cite{MiSu} that meet in one point, the cuspidal point of the
quintic.

These two curves generate a pencil $\cL_{\PP^2}=\langle
5C_1,2C_2\rangle$ of rational curves of degree $10$ with a unique
base point such that $\PP^2\setminus (C_1\cup C_2)$ is a cylinder.
Similarly as in the proof above, every del Pezzo surface $Y$ of
degree $d\ge 4$ can be obtained, along with a $(-K_Y)$-polar
cylinder, by blowing up a certain set of $9-d$ points on $C_1$.
Indeed, we can write
$$-K_{\PP^2}\equiv (\textstyle{\frac{3}{2}}-\varepsilon) C_1
+\frac{2}{5}\varepsilon C_2$$ with an appropriate $\varepsilon>0$,
and then proceed in the same fashion as in the proof. \eexa

\bexa \label{ex10} Picking up four points $P_1,\ldots,P_4$ in
$\PP^2$ in general position, we consider the pencil of lines
$\cL_{\PP^2}$ on $\PP^2$ generated by $l_1=(P_1P_2)$ and
$l_2=(P_3P_4)$. The blowup $\sigma: Y\to \PP^2$ of these points
yields a del Pezzo surface $Y$ of degree $5$. We have
$$-K_{\PP^2}\equiv\frac{3}{2}l_1+\frac{3}{2}l_2\qquad\mbox{ and
so}\qquad
-K_Y\equiv\frac{3}{2}l'_1+\frac{3}{2}l'_2+\frac{1}{2}\sum_{i=1}^4
E_i,$$ where $l_i'$ is the proper transform of $l_i$, $i=1,2$, and
$E_i$ is the exceptional $(-1)$-curve over $P_i$, $i=1,\ldots,4$.
Then $L^{(1)}=l_1'+E_1+E_2$ and $L^{(2)}=l_2'+E_3+E_4$ are the
only degenerate fibers of the pencil $\cL=\sigma^{-1}_*
(\cL_{\PP^2})$ on $Y$. Since
$$D:=\frac{1}{2}(L^{(1)}+L^{(2)})+(l_1'+l_2')\equiv-K_Y\,,$$ the open set
$$Y\setminus \supp (D)=\PP^2\setminus (l_1\cup l_2)\simeq \A^1_*\times\A^1$$
is a $(-K_Y)$-polar cylinder on $Y$. A similar construction can be
applied to any del Pezzo surface of degree $d\ge 5$. \eexa

\bexa\label{ex50} Consider the pencil $\cL_{\PP^2}$ of unicuspidal
rational curves $\alpha yz^{n-1}+\beta x^n=0$ in $\PP^2$, where
$n\ge 1$. Blowing up $k\le 4$ points in $\PP^2$, at most two on
each of the lines $x=0$ and $y=0$ off their common point $(0:0:1)$
we obtain examples of  $(-K_Y)$-polar cylinders on an arbitrary
del Pezzo surface of degree $d\ge 5$. For $n=1$ and $k=4$ this
gives again the cylinder from Example \ref{ex10}. \eexa

\brem\label{not3} The idea to start with a $(-K_{\PP^2})$-polar
cylinder on $\PP^2$ cannot be carried out any more in case of a
smooth cubic surface $Y\subseteq\PP^3$. Indeed, suppose we are
given a cylinder $\PP^2\setminus C\simeq Z\times\A^1$, where $C$
is a reduced plane curve of degree $d$, not necessarily smooth or
irreducible, and let $\cL_{\PP^2}$ be the corresponding pencil.
Then $\Bs(\cL_{\PP^2})$ consists of one point $P_0\in C$, and
$C\setminus\{P_0\}$ is a disjoint union of components isomorphic
to $\A^1$. Performing a blowup $\sigma:Y\to\PP^2$ of $m$ points
$P_i\in C\setminus\{P_0\}$ with exceptional curves
$E_j=\sigma^{-1}(p_j)$, $i=1,\ldots,m$, from the equalities
$$3\sigma^*(C)\thicksim d\sigma^* (-K_{P^2})=-dK_Y+d\sum_{j=1}^m
E_j\quad\mbox{and}\quad 3\sigma^*(C)=3C'+3\sum_{j=1}^m E_j\,,$$
where $C'$ is the proper transform of $C$ on $Y$, we obtain
$$-dK_Y\thicksim 3C'+(3-d)\sum_{j=1}^m E_j=:D\,.$$
Here $D$ is an effective divisor with $\supp (D)=C'+\sum_{j=1}^m
E_j$ if and only if $d\le 2$ i.e., $C$ is a line or a conic. Since
the centers of blowup $P_i$, $i=1,\ldots,m$, are situated on $C$
and $Y$ must be del Pezzo, we have $m\le 5$ and so $\deg (Y)\ge
4$. \erem

In the next example, starting with a pencil on $\PP^2$ with five
base points, we construct  a $(-K_Y)$-polar cylinder of different
type on arbitrary del Pezzo surface $Y$ of degree $d=5$,
 by
resolving all the base points but one.

\bexa\label{ex11} Consider the following pencil $\cL_{\PP^2}$ of
rational plane sextics:
\[
\alpha (y^2z-x^3)^2+\beta (y^2-xz)(y^4-x^4)=0\,.
\]
The base locus of $\cL_{\PP^2}$ consists of the points $P_0,\dots,
P_4$, where $P_0=(0:0:1)$ and $\{P_1,\dots, P_4\} =(x^2=z^2,\quad
xz=y^2)$. Furthermore, $\cL_{\PP^2}$ has no fixed component and so
its general member $L$ is irreducible. Since $\mult_{P_{0}} (L)=4$
and $\mult_{P_{i}} (L)=2$ for $i=1,\dots, 4$, the curve $L$ is
rational. Any singular point $P_i$ of $L$ is resolved by one
blowup, and the singularity of $L$ at $P_0$ is cuspidal. No three
of the points $P_1,\dots, P_4$ are collinear. Therefore the blowup
$\sigma : Y\to \PP^2$ of the latter points yields a del Pezzo
surface $Y$ of degree $5$, and any such surface arises in this
way. Let $\cL$ be the proper transform of $\cL_{\PP^2}$ on $Y$,
and let $P=\sigma^{-1}(P_0)$. Then $\cL$ is a pencil of rational
curves with a cuspidal singularity at the unique base point $P$,
smooth and disjoint outside $P$.

There are exactly two degenerate members of $\cL_{\PP^2}$, namely
the double cuspidal cubic $C'=2(y^2z=x^3)$ and the union $C''$ of
the conic $(y^2=xz)$ and the four lines $(y^4=x^4)$. We have
$-K_{\PP^2}\equiv D_{\PP^2}:=\frac14 C'+\frac14 C''$. Let $D$ be
the proper transform of $D_{\PP^2}$ on $Y$. Since $\mult_{P_{i}}
(D_{\PP^2}) =1$ for $i=1,\dots, 4$, we have
$K_Y+D=\sigma^*(K_{\PP^2}+D_{\PP^2})\equiv 0$. Hence $U=Y\setminus
(C'\cup C'')$  is  a $(-K_Y)$-polar cylinder on $Y$ (see also
Example \ref{ex111} below).\eexa

\subsection{The Makar-Limanov invariant on affine cones over del Pezzo surfaces}
\bsit\label{mlinv}For an algebra $A$ over a field $k$, its
Makar-Limanov invariant $\ML(A)$ is defined as the intersection of
the kernels of all locally nilpotent derivations on $A$. It is
trivial if $\ML(A)=k$. Following \cite{MiMa} we say that $A$ is of
class $\ML_i$ if the quotient field $\Frac\, (\ML(A))$ has
transcendence degree $i$. If $\ML(A)$ is finitely generated then
$i=\dim (Z)$, where $Z=\spec \ML(A)$. Thus $A\in\ML_0$ whenever
$A$ has trivial Makar-Limanov invariant. For instance, $\A^3\in
\ML_0$ (regarded as the affine cone over $\PP^2$).

For $A$ graded there are graded versions $\ML^{(h)} (A)$ and
$\ML^{(h)}_i$ of $\ML (A)$ and $\ML_i$, respectively \cite{FZ2},
where one restricts to homogeneous locally nilpotent derivations.
Clearly, $\ML (A)\subseteq \ML^{(h)} (A)$. Hence the usual ML
invariant is trivial if the homogeneous is. \esit

\bthm\label{mli} Let $X$ be the affine cone  over a smooth,
anticanonically embedded del Pezzo surface $Y\subseteq\PP^d$ of
degree $d\ge 4$. Then the homogeneous Makar-Limanov invariant
$\ML^{(h)}(X)$ is trivial i.e., $X\in \ML^{(h)}_0$. \ethm

\bproof  By Theorem \ref{dpaut} and Proposition \ref{akankl}, for
$d\ge 6$ the surface $Y$ and the cone $X=\affcone_{-K_Y} (Y)$ are
toric. Since $X$ is not isomorphic to a product $X'\times\A^1_*$,
by Lemma 4.5 in \cite{Li} the homogeneous Makar-Limanov invariant
of $X$ is trivial.

It remains to show that $X\in \ML^{(h)}_0$ for $d=4,5$. Note that
for an arbitrary graded algebra $A=\bigoplus_i A_i$, the graded
subalgebra $\ML^{(h)}(A)$ is non-trivial if and only if there
exists a non-constant homogeneous element $h\in A_n\cap
\ML^{(h)}(A)$ (so $h$ is annihilated by all homogeneous locally
nilpotent derivations on $A$). In the case of an affine cone $X$,
the degree $n=\deg (h)$ is positive. Hence $\Gamma=\V(h)\in
|n(-K_Y)|$ is an effective ample divisor on $Y$.

Let as before $Y\subseteq\PP^d$ be a del Pezzo surface of degree
$d\ge 4$. Suppose on the contrary that $\ML^{(h)}(X)\neq \C$, and
let $h\in A_n\cap \ML^{(h)}(A)$ be nonconstant. Then the affine
cone over the curve $\supp\,(\Gamma)$ is a divisor on $X$ stable
under any $\C_+$-action defined by a homogeneous locally nilpotent
derivation on $\cO(X)$. Hence for every $(-K_Y)$-polar cylinder
$U$ on $Y$, the curve $\supp\,(\Gamma)$ consists of components of
the members of the  linear pencil $\cL$ on $Y$ associated with
$U$. In particular, for all $(-K_Y)$-polar cylinders on $Y$ there
must be a common component of the associated linear pencils $\cL$.

For $d=5$, we let $\sigma:Y\to\PP^2$ be the blowup of four points
$P_1,\ldots,P_4$ in $\PP^2$ with exceptional curves
$E_i=\sigma^{-1}(P_i)$. There are exactly ten lines on $Y$.
Besides $E_1,\ldots, E_4$ these are the proper transforms $l_{ij}$
of the lines $(P_iP_j)$ on $\PP^2$, where $1\le i<j\le 4$. For
every pair of lines $(l_{ij},l_{i'j'})$ with distinct indices
$i,j,i',j'$, the curves
$$L_1:=l_{ij}+E_i+E_j\quad\mbox{and}\quad
L_2:=l_{i'j'}+E_{i'}+E_{j'}$$ are the only degenerate members of a
cylindrical linear pencil on $Y$ (cf.\ Example \ref{ex10}). The 3
such pencils have no common component except for the lines
$E_1,\ldots, E_4$.

Let us replace the lines $E_1,\ldots,E_4$ on $Y$ by some other
four disjoint lines, e.g. by $l_{12},l_{13},l_{23},E_4$. We
consider also the three associated cylindrical pencils on $Y$
e.g., that with degenerate members
$$L'_1:=E_2+l_{12}+l_{23}\quad\mbox{and}\quad
L'_2:=E_4+l_{13}+l_{24}\,.$$ These pencils have no common
component except for
 the lines $l_{12},l_{13},l_{23},E_4$.
 The line $E_4$ is the only common component
of all six above pencils. With yet further choice of a pencil, we
can eliminate also this latter line. Thus the homogeneous
Makar-Limanov invariant of $Y$ is trivial, as stated.

Let further $d=4$, and let $\sigma_0:Y\to\PP^2$ be the blowup of
five points $P_1,\ldots,P_5$ in general position in $\PP^2$, with
exceptional curves $E_i=\sigma^{-1}(P_i)$,  $i=1,\ldots,5$. We let
$C$ denote the unique smooth conic through the points $P_i$. Given
a point $Q\in C$ different from the $P_i$, similarly as in the
proof of Theorem \ref{cuco} we consider the pencil of conics on
$\PP^2$ generated by $C$ and $2l_Q$, where $l_Q$ is the tangent
line to $C$ at $Q$. Two different such pencils on $\PP^2$ have no
common member except for the conic $C$ itself. Thus $C'$ and the
lines $E_i$, $i=1,\ldots,5$, are the only common components of the
induced cylindrical pencils on the del Pezzo surface $Y$ of degree
4.

Consider next the contraction $\sigma_1:Y\to \PP^2$ of the five
disjoint lines $C',l_{12},l_{13},l_{14},l_{15}$
 on $Y$. Then $\sigma_1(E_1)$ is a conic in $\PP^2$,
 which plays now the
role of $C$. Once again, the only common components of the induced
cylindrical pencils on $Y$ are $E_1$ and the five disjoint lines
above meeting $E_1$.

Likewise, for the six different contractions $\sigma_i:Y\to
\PP^2$, $i=0,\ldots,5$, the only common component of the induced
cylindrical pencils on $Y$ is $C'$. However, the ample divisor
$\Gamma$ as at the beginning of the proof cannot be supported by
$C'$. This contradiction finishes the proof. \eproof

\noindent {\bf Problem.} {\em Describe all affine cones whose
homogeneous Makar-Limanov invariant is trivial.}

\section{On existence of $\C_+$-actions on cones over cubic surfaces}
In this section we analyze in detail the case of a smooth cubic
surface $Y\subseteq \PP^3$. We do not know whether the affine cone
$X=\affcone (Y)$ carries a $\C_+$-action. However, we obtain in
Proposition \ref{sumup} and Theorem \ref{nefna} below a detailed
information on an eventual anticanonical polar cylinder in $Y$.
This makes the criterion \ref{ext2} of the existence of a
$\C_+$-action much more concrete in our particular case. We adopt
the following convention.

\bsit\label{1000} {\bf Convention.} Let $Y$ be a smooth cubic
surface in $\PP^3$. Suppose that the affine cone $X=\affcone
(Y)\subseteq\A^{4}$ admits an effective $\C_+$-action. Then by
Corollary \ref{ext2} $Y$ possesses a $(-K_Y)$-polar cylinder
$U\simeq \A^1\times Z$, where $Z$ is an affine smooth rational
curve. In other words $Y\setminus U=\supp (D)$, where\be\label{00}
D=\sum_{i=1}^n\delta_i\Delta_i\equiv -K_Y \ee is an effective
$\Q$-divisor on $Y$ with $\delta_i\in\Q_{> 0} \,\forall
i=1,\ldots,n$.

\subsection{The linear pencil on a cubic surface compatible with a cylinder}
Let $\cL$ be the pencil on $Y$ with general member
$L_z=\overline{{\rm pr}_2^{-1}(z)}$ for $z\in Z$. It is easily
seen that $\cL$ has at least one degenerate member. In what
follows we suppose that $\supp D$ does not contain a
non-degenerate member of $\cL$ (otherwise, up to numerical
equivalence, we replace such a member by a degenerate one). Under
these assumptions, the following hold. \esit

\blem\label{picard} The support of $D$ is connected and simply
connected, and contains at least $7$ irreducible components.\elem
\bproof The projection $\pr_2:U\to Z$ extends to a rational map
$Y\dashrightarrow \PP^1$ defined by the pencil $\cL$ as above. A
general member $L$ of $\cL$ is a rational curve smooth off a
unique point $P$, where $\{P\}=L\cap \supp (D)$, and
$L\setminus\{P\}=L\cap U\simeq\A^1$. Thus $L$ is unibranch at $P$.
The only base point of $\cL$ (if exists) is contained in $\supp
(D)$.

Since $D$ is ample, by the Lefschetz Hyperplane Section Theorem
$\supp (D)$ is connected. Resolving, if necessary, the base point
of $\cL$ by a modification $p:W\to Y$ yields a rational surface
$W$ with a pencil of rational curves $\cL_W={p_*}^{-1} \cL$ that
are fibers of $q=q_{|\cL_W|}:W\to \PP^1$. By Zariski's Main
Theorem, the total transform $p^{-1}(\supp D)$ is still connected,
and is a union of a $(-1)$-section, say, $S$ of $\cL_W$ and of
some number of rational trees contained in fibers of $\cL_W$.
Hence $p^{-1}(\supp D)$ is also a tree of rational curves i.e.,
is connected and simply connected. The exceptional divisor $E$ of
$p$ being a subtree of $p^{-1}(\supp D)$, the contraction of $E$
does not affect the simply-connectedness. This proves the first
assertion.

Since $Z$ is a rational smooth affine curve, we have $\Pic
(U)=\Pic (\A^1)\times\Pic (Z)=0$ . By virtue of the exact sequence
\be\label{es}  G:=\sum_{i=1}^n\Z\Delta_i\to\Pic (Y)\to\Pic
(U)=0\,\ee the free abelian group $G$ generated by the components
$\Delta_i$ of $D$ surjects onto  $\Pic (Y) \simeq \Z^7$. Therefore
$ {\rm rk} (G)\ge 7$, which proves the second assertion. \eproof

\blem\label{0} The pencil $\cL$ has a unique base point, say, $P$,
and $\deg (\cL)\ge 3$. \elem

\bproof If on the contrary $\Bs(\cL)=\emptyset$, then the pencil
of conics $\cL$ on $Y$ with a section, say, $S=\Delta_0$ defines a
morphism $\varphi_{|C|}:Y\to\PP^1$ (extending the projection
$\pr_2$ of the cylinder) with exactly five degenerate fibers
$L_1,\ldots,L_5$. Each
 degenerate fiber consists of two lines on $Y$ intersecting
transversally at one point. At most one of these two lines, say,
$l_i$ meets the cylinder $U$, while the other one, say, $\Delta_i$
is a components of $D$. Since $D$ is connected we have
$\Delta_i\cdot S=1$ and $l_i\cdot S =0$, $i=1,\ldots,5$. By the
Adjunction Formula we get
$$1=(-K_Y)\cdot l_i=D\cdot l_i=\delta_i-x,$$
where $x=0$ if $l_i\neq\Delta_j$ $\forall j$ and $x=\delta_j>0$
otherwise. Hence $\delta_i= 1+x\ge 1$. Similarly, for a general
fiber $L$ of $\cL$,
$$2=(-K_Y)\cdot L=D\cdot L =\delta_0S \cdot L =\delta_0$$
and so $\delta_0=2$.

On the other hand, by the Adjunction Formula
$$-K_Y\cdot S =2+S ^2=D\cdot S = 2S ^2+\sum_{i=1}^5 \delta_i\,.$$ Therefore,
$$2=S ^2+\sum_{i=1}^5 \delta_i=S ^2+\sum_{i=1}^5 \delta_i\ge -1+5=4,$$
a contradiction. The inequality $(-K_Y)\cdot\cL\ge 3$ is now
immediate.\eproof

\brems\label{345}1. Actually the degree of $\cL$ must be
essentially higher, since by Lemma \ref{070} below $D$ has $8$
irreducible components.

2. The assertion of the lemma holds also for any del Pezzo surface
of degree 4 or 5, with a similar proof. However, it fails for
degree 6. Indeed, pick 3 points $P_0$, $P_1$, $P_2$ in general
position in $\PP^2$, and consider the pencil generated by the
lines $l_i=(P_0P_i)$, $i=1,2$. Blowing up these points we get a
del Pezzo surface $Y$ of degree 6 with a base point free pencil.
Then the complement in $Y$ of the total transform of $l_1\cup l_2$
is a $(-K_Y)$-polar cylinder.\erems

\bsit\label{000} In the sequel we frequently use the following
commutative diagram:

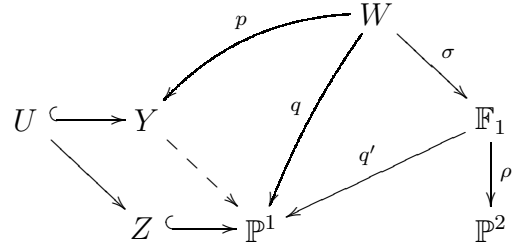
\be\label{di} \xymatrix{
&&&&W\ar@/_1pc/[dll]_{p}\ar[dr]^{\sigma}\ar@/_0.3pc/[ddl]_q&
\\
&U\ \ar@{^{(}->}[r]\ar[dr]&
Y\ar@{-->}[dr]&&&\F_1\ar[dll]_{q'}\ar[d]^{\rho}
\\
&&Z\ \ar@{^{(}->}[r]&\PP^1&&\PP^2& } \ee\medskip

\noindent where $p:W\to Y$ is the minimal resolution of the base
point $P$ of $\cL$, $q:W\to\PP^1$ is the induced pencil,
$\sigma:W\to \F_1$ is composed of the contraction of all
components of degenerate fibers of $q$ except those which meet the
exceptional $(-1)$-section $S_W$ of $q$, and $\rho:\F_1\to\PP^2$
is the contraction of this exceptional section. \esit

\blem\label{delta-i} $\delta_i< 1$ in \eqref{00} for all
$i=1,\ldots,n\,.$\elem\bproof It is enough to show that $\delta_1<
1$. By symmetry then $\delta_i< 1$ $\forall i$. Since the
anticanonical divisor of $Y$ is ample, we have $(-K_Y)\cdot
\Delta_i
> 0$ $\forall i=1,\ldots,n$. We distinguish between the following
3 cases : \bnum\item[(i)] $(-K_Y)\cdot \Delta_1\ge 3$,
\item[(ii)] $(-K_Y)\cdot \Delta_1=2$, and \item[(iii)] $(-K_Y)\cdot
\Delta_1=1$. \enum
In case (i) suppose on the contrary that $\delta_1\ge 1$. Since
$n\ge 7$  by Lemma \ref{picard} and the divisor  $-K_Y\equiv D$ is
ample, we obtain
$$3=(-K_Y)\cdot D=\sum_{i=1}^n\delta_i (-K_Y)\cdot \Delta_i>\delta_1 (-K_Y)\cdot
\Delta_1\ge 3,$$ which gives a contradiction.

In case (ii) $\Delta_1\subseteq Y$ is a conic, $\Delta_1^2=0$, and
$-K_Y\equiv\Delta_1+E$, where $E$ is the residual line cut out on
$Y$ by a plane in $\PP^3$ through $\Delta_1$.

Let $E=\Delta_i$ for some $i>1$; we may assume that $i=2$. Then
$-K_Y\equiv \Delta_1+\Delta_2$, $\Delta_1^2=0$ and
$\Delta_1\cdot\Delta_2=2$. Thus \be\label{eqa1}
2=(-K_Y)\cdot\Delta_1=D\cdot\Delta_1\ge
\delta_2\Delta_1\cdot\Delta_2 =2\delta_2,\ee so $\delta_2\le 1$.
Furthermore,
$$1=\Delta_2\cdot D=\sum_{i=1}^n \delta_i\Delta_2\cdot\Delta_i\ge 2\delta_1-
\delta_2\,.$$ Hence $\delta_1\le\frac{1}{2} (\delta_2+1)\le 1$.

If $\delta_1=1$ then also $\delta_2=1$. Since
$\Delta_1+\Delta_2\equiv-K_Y\equiv D$ we get
$D=\Delta_1+\Delta_2$. Thus $n=2$, which contradicts Lemma
\ref{picard}. Therefore in this case $\delta_1<1$, as stated.

If further $E\neq\Delta_i$ $\forall i$ then \be\label{eqa2}
1=(-K_Y)\cdot E=D\cdot E\ge \delta_1\Delta_1\cdot E=2\delta_1, \ee
hence $\delta_1\le 1/2$. Thus anyway, $\delta_1<1$ in case (ii).

In case (iii) $\Delta_1$ is a line on $Y$. Let $C$ be a residual
conic of $\Delta_1$, so that $\Delta_1+C\equiv-K_Y$ is a
hyperplane section. We have as before $$2=(-K_Y)\cdot C=D\cdot
C\ge \delta_1\Delta_1\cdot C=2\delta_1,$$ hence $\delta_1\le 1$.

If $\delta_1=1$ then $D\cdot C=\Delta_1\cdot C=2$ and so
$(D-\Delta_1)\cdot C=0$. Therefore the divisor $\supp
(D-\Delta_1)$ is supported on the members of the pencil of conics
$|C|$ on $Y$. The curve $\Delta_1$ meets each fiber twice, and so
the morphism $\varphi_{|C|}$ restricted to $\Delta_1$ has 2 branch
points.

By Lemma \ref{picard} the curve $\supp D$ is simply connected,
hence it cannot contain the whole fiber of $\varphi_{|C|}$ which
meets the component $\Delta_1$ of $D$ at two distinct points. We
claim however that if a degenerate fiber $l_1+l_2$ of the pencil
$|C|$ contains a component, say, $\Delta_i=l_1$ of $D$, then its
second component $l_2=\Delta_j$ is also contained in $\supp D$
and, moreover, $\delta_i=\delta_j$. Indeed, since $\delta_1=1$ and
$\Delta_1$ is a line on $Y$ we have \be\label{010} 1=\Delta_i\cdot
D=\Delta_i\cdot
 \Delta_1+ \Delta_i\cdot (D-\Delta_1)=1-\delta_i+\sum_{k\neq 1,i}\delta_k
 \Delta_i\cdot \Delta_k\,.\ee
 The only component of the latter sum that meets $\Delta_i$
 can be the line $l_2$. Hence $l_2=\Delta_j$ for some $j\neq 1,i$.
 Now \eqref{010} shows that $\delta_i=\delta_j$, as claimed.

Furthermore, since $\Delta_i\cup \Delta_j$ meets $\Delta_1$ twice
and $\supp D$ is a tree, the line $\Delta_1$ passes through the
intersection point $\Delta_i\cap \Delta_j$. On the other hand,
$\Delta_1$ is tangent to exactly two members of the pencil $|C|$,
which are
 either smooth or consist of two lines
$\Delta_i$ and $\Delta_j$ meeting $\Delta_1$ at their common point
(an Eckardt point of $Y$). By the simply connectedness of $\supp
D$, none of the other components of members of $|C|$ can be
contained in $\supp D$. Hence $\supp D$ can contain at most 5
components, namely, $\Delta_1$ and the components of two
degenerate members tangent to $\Delta_1$. However, this
contradicts Lemma \ref{picard}, since by this lemma $\supp D$
consists of at least $7$ components. Now the proof is completed.
\eproof

\blem\label{pass} Every component of the degenerate members of the
pencil $\cL$ on $Y$ passes through the base point $P$ of
$\cL$.\elem

\begin{proof}
Assume on the contrary that there is a component $C_0$ of a
degenerate member $L^{(0)}$ of $\cL$ such that $P\notin C$. By
Zariski's Lemma $C_0^2<0$. Hence $C_0$ is a $(-1)$-curve on $Y$
and so $D\cdot C_0=(-K_Y)\cdot C_0=1$. By an easy argument (cf.\
the proof of Lemma \ref{picard}) the curve $\supp L^{(0)}$ is
connected and simply connected (this is a tree of rational curves
outside $P$). If there is another component $C_1$ of $L^{(0)}$
which meets $C_0$ and does not pass through $P$, then $L^{(0)}$
contains the configuration of two crossing lines $C_0+C_1$ which
do not pass through $P$. Then $\cL$ must be the linear system of
conics $|C_0+C_1|$. By Lemma \ref{0} this leads to a
contradiction. Hence $C_0$ cannot separate $\supp L^{(0)}$ and so
$C_0$ meets the complement $\supp ({L^{(0)}- C_0})$ at one point
transversally.

It follows that $\Delta_i\cdot C_0=1$ for a unique index $i$. If
$C_0$ is not a component of $D$, then $1=D\cdot C_0= \delta _i$,
which contradicts Lemma \ref{delta-i}. Similarly, if
$C_0=\Delta_j$, then $1=D\cdot C_0= \delta _i-\delta_j$. Hence
$\delta_i=1+\delta_j>1$. Once again, this contradicts Lemma
\ref{delta-i}.
\end{proof}

\blem\label{050} The pencil $\cL$ has at most two degenerate
members. \elem

\bproof Recall (see \ref{000}) that $p: W\to Y$ stands for the
minimal resolution of the base locus of $\cL$ and $q: W\to \PP^1$
for the fibration given by $p^{-1}_*\cL$. Write $p$ as a
composition of blowups of points over $P$: \be\label{bu} p:
W\stackrel{p_1}{\longrightarrow}
W_1\stackrel{p_2}{\longrightarrow} \cdots
\stackrel{p_N}{\longrightarrow} W_N=Y, \ee where the exceptional
divisor $S_W$ of $p_1$ is a $q$-horizontal $(-1)$-curve on $W$. A
general fiber $L$ of $q$ is a smooth rational curve meeting $S_W$
at one point. Indeed, $L\setminus S_W\simeq p(L)\setminus P\simeq
\A^1$. Therefore $S_W$ is a section of $q$.

Let $C_1,\dots, C_m$ be the components of degenerate fibers
$F_1,\dots, F_m$ of $q$ meeting $S_W$. We claim that all the
curves $C_i$ are $p$-exceptional. Indeed, otherwise for some $i$,
the image $p(\overline{(F_i\setminus C_i)})$ would be a component
of a degenerate member of $\cL$ which does not pass through $P$.
The latter contradicts Lemma \ref{pass}.

Note that, on each step, the exceptional divisor of $p_k\circ
\ldots \circ p_N$ is an SNC tree of rational curves. On the other
hand, all the curves $p_1(C_i)$ on $W_1$ pass through the point
$p_1(S_W)$. Therefore $m\le 2$. \eproof

\blem\label{070} The pencil $\cL$ has exactly two degenerate
members, say, $L^{(1)}$ and $L^{(2)}$. Furthermore, $\supp D=\supp
(L^{(1)}+L^{(2)})$ consists of $8$ irreducible components i.e.,
$n=8$ in \eqref{00}. \elem \bproof Assume that the only degenerate
member of $\cL$ is $L^{(1)}$. In this case, $\supp D\subseteq
\supp L^{(1)}$, $Z\supseteq\A^1$, and $U\supseteq \A^2$.

If $Z\simeq\A^1$ then $\Pic (U)=0$ and $H^0(U,\cO_U)=\C$, hence
$\Pic (Y)\simeq \sum_{i=1}^n \Z\cdot \Delta_i$. Thus in this case
$n=7$ and $-K_Y= \sum_{i=1}^7 m_i\Delta_i$ for some $m_i\in \Z$.
On the other hand, $-K_Y\equiv D=\sum_{i=1}^7 \delta_i\Delta_i$,
where $0<\delta_i<1$ $\forall i$ according to Lemma \ref{delta-i}.
Since the decomposition of $-K_Y$ in $\Pic (Y)\otimes\Q$ is
unique, this yields a contradiction.

If further $Z\simeq\PP^1$ then $\Pic (U)\simeq\Z$ and
$H^0(U,\cO_U)=\C$. From the exact sequence \eqref{es}, where $\Pic
(U)=0$ is replaced by $\Pic (U)\simeq\Z$, we obtain $n=6$.  By
Lemma \ref{picard}, this leads again to a contradiction. Therefore
$\cL$ has indeed two degenerate members.

As for the second assertion, assuming on the contrary that $\supp
D\neq \supp (L^{(1)}+L^{(2)})$ we would have $Z\supseteq\A^1$. Now
the same argument as before yields a contradiction. Since the
Picard group $\Pic (Y)\cong\Z^7$ is generated by the irreducible
components of $L^{(1)}+L^{(2)}$ and $L^{(1)}\equiv L^{(2)}$ is the
only relation between these components, we obtain that $n=8$.
\eproof

\bcor\label{100} The pencil $\cL$ is ample. Furthermore, for every
irreducible component $C$ of a member of the pencil $\cL$ on $Y$
we have $C\setminus\{P\}\simeq\A^1$, and two such components have
just the point $P$ in common. \ecor \bproof The first assertion
follows immediately from Lemma \ref{pass} by the Nakai-Moishezon
criterion. As for the second one, it follows from the well known
fact that on an affine surface $V=Y\setminus L$, where $L$ is a
general fiber of $\cL$, every degenerate fiber of the
$\A^1$-fibration $\varphi_{|\cL|}:V\to\A^1$ is a disjoint union of
affine lines, see e.g., \cite{Mi}, \cite{Za}. \eproof

\blem\label{020} The pair $(Y,D)$ is not log canonical at $P$.
\elem \bproof Let $D_W$ denote the crepant pull-back of $D$ on $W$
 as in (\ref{di}) i.e., a $\Q$-divisor on $W$ such that
\[
K_W+D_W=p^*(K_Y+D)\quad\mbox{and} \quad p_*D_W=D\,.
\]
The exceptional $(-1)$-section $S$ on $W$ is the only
$q$-horizontal component of $D_W$. For a general fiber $L$ of $q$
we have $2=(-K_W)\cdot L=D_W\cdot L$. Therefore, the discrepancy
$a(S,D)$ (i.e., the coefficient of $S$ in $D_W$ with the opposite
sign) equals $-2$. This proves our assertion. \eproof

\bcor\label{030} $\mult_P (D)>1$. \ecor \bproof If $\mult_P (D)\le
1$, then the pair $(Y,D)$ is canonical at $P$, because $P\in Y$ is
a smooth point. In particular, it is log canonical. This
contradicts Lemma \ref{020}. \eproof

\blem\label{040} Any line $l$ on $ Y$ through $P$ is contained in
$\supp D$. \elem \bproof Assuming the contrary we obtain
$1=(-K_Y)\cdot l=D\cdot l\ge \mult_P (D) >1$, a contradiction.
\eproof

\blem\label{080} $P$ cannot be an Eckardt point on $Y$. \elem
\bproof Suppose the contrary. Then by Lemma \ref{040}, up to a
permutation we may assume that $\Delta_1$, $\Delta_2$, $\Delta_3$
are lines through $P$, where $\delta_1\le\delta_2\le\delta_3$.
Since $D\cdot\Delta_1=(-K_Y)\cdot\Delta_1=1$, for an effective
$\Q$-divisor on $Y$ \be\label{dva}
D':=\frac{1}{1-\delta_1}(D-\delta_1(\Delta_1+\Delta_2+
\Delta_3))\,\ee we obtain $D'\cdot\Delta_1=1$  and by (\ref{dva})
$$D\equiv-K_Y\equiv \Delta_1+\Delta_2+ \Delta_3\equiv D'\,.$$ Now the
proof of Lemma \ref{020} works equally for $D'$. Hence the pair
$(Y,D')$ is not canonical at $P$ and so $\mult_P (D')>1$. This
contradicts the inequality $\mult_P (D')\le D'\cdot\Delta_1=1$.
\eproof

\blem\label{060} The fibration $q:W\to\PP^1$ has exactly two
degenerate fibers. The general member $L$ of $\cL$ is singular at
$P$.\elem

\bproof The first assertion follows from Lemma \ref{070}. Let us
show the second one. Assuming the contrary, for a smooth rational
curve $L$ on $Y$ we have by adjunction $(K_Y+L)\cdot L=-2$. The
Mori cone of $Y$ being spanned by the $(-1)$-curves
$E_1,\ldots,E_{27}$, there is a decomposition
$$L\equiv \sum_{i=1}^{27} \alpha_i E_i,$$
where $\alpha_i \in\Q_{\ge 0}$. Hence $(K_Y+L)\cdot E_i<0$ for
some $i$. Thus $L\cdot E_i<(-K_Y)\cdot E_i=1$, and so $L$ cannot
be ample. This contradicts Corollary \ref{100}. \eproof

\brem \label{w} The minimal resolution $p:W\to Y$ of the base
point $P$ of $\cL$ dominates the embedded minimal resolution of
the cusp at $P$ of a general member $L$ of $\cL$. The exceptional
divisor $E=p^{-1}(P)\subseteq W$ is a rational comb with the
number of teeth equal the length of the Puiseaux sequence of the
cusp. The only $(-1)$-curve $S_W$ in $E$ is sitting on the handle
of the comb. Hence $E=E^{(1)}+S_W+E^{(2)}$, where
$E^{(k)}\subseteq p^{-1}(L^{(k)})$, $k=1,2$, and exactly one of
the $E^{(k)}$ is a negatively definite linear chain of rational
curves.

The degenerate members $L_W^{(1)}$, $L_W^{(2)}$ of the induced
linear system $\cL_W$ on $W$ have the following structure:
$L_W^{(k)}$ consists of $E^{(k)}$ and the proper transforms
$\Delta_i'$ of the components $\Delta_i\subseteq L^{(k)}$, called
{\em feathers}. The feathers are disjoint in $W$, and each of them
meets $E^{(k)}$ at one point transversally. Let us illustrate this
for the surface $Y$ from Example \ref{ex11}, with the pullback
$\cL_W$ of the pencil $\cL$ on $Y$.\erem

\bexa\label{ex111} Up to interchanging the degenerate fibers
$L_W^{(1)}$ and $L_W^{(2)}$ of $\cL_W$, their structure is as
follows. The fiber $L_W^{(1)}$ contains a component $E_1$ with
self-intersection $-5$ joint with the $(-1)$-section $S_W$, and
also the $(-1)$-feathers $\Delta'_1,\ldots,\Delta_5'$ meeting
$E_1$. The second fiber $L_W^{(2)}$ contains the only feather
$\Delta'_0$ of multiplicity 2. The weighted dual graph of the
configuration $L_W^{(1)}+S_W+L_W^{(2)}$ is as follows:

\bigskip

\be\label{graph0}\; \cou{\!\!\!\!\!\!\!\!\!\!\!\!
-5}{E_1}\nlin\boxshiftup{\!\!\!\!
(-1)_5}{\,\,\,\,\,\,\,\,\,\,\,\,\,\,\,\,\,\,\,\,(\Delta_i')_{i=1,...,5}}
\llin\cou{-1}{S_W}\llin \cou{-2}{E_2}\llin
\cou{\!\!\!\!\!\!\!\!\!\!\!\! -2}{E_3} \nlin\cshiftup{\!\!\!\!
-1}{ \Delta'_0}\llin\cou{-2}{E_{4}}\ee

\medskip

\noindent where the box denotes a disjoint union of five
$(-1)$-feathers $\Delta'_1,\ldots,\Delta'_5$ joint with $E_1$. The
exceptional divisor of $p: W\to Y$ is $E=E_1+\ldots +E_4$. While
$\Delta_i=p(\Delta'_i)$ ($i=0,\ldots,5$) are the components of the
degenerate fibers $L^{(j)}=p(L_W^{(j)})$ ($j=1,2$) of the pencil
$\cL$ on $Y$. More precisely, $\Delta_0$ is the proper transform
in $Y$ of the cuspidal cubic $C'$, $\Delta_1$ is that of the conic
$C''$, while $\Delta_2,\ldots,\Delta_5\subseteq Y$ are the proper
transforms of the four lines $(y^4=x^4)$ in $\PP^2$. \eexa

\bsit\label{lisi} Let $Q=(\rho\circ \sigma)(S_W)\in\PP^2$, where
as in \ref{000} $S_W$ stands for a $(-1)$-section of $q$ contained
in $\Exc(p)$. Then $\cL_{\PP^2}:=\rho_* (\sigma_*\cL_{W})$ is the
linear pencil of lines through $Q$ on $\PP^2$. We let
$\varphi:=\rho\circ\sigma\circ p^{-1}: Y\dashrightarrow \PP^2$,
and we let $\cH$ denote the proper transform of
$\cH_{\PP^2}:=|\cO_{\PP^2}(1)|$ on $Y$ via $\varphi$. With this
notation, the following holds. \esit

\blem\label{lici1} $\cL\subseteq \cH$. \elem \bproof We have
\[
\cH_W:=(\rho \circ \sigma)^{-1}_*\cH_{\PP^2}= (\rho\circ\sigma)^*
\cH_{\PP^2}= \sigma^*(\rho^* \cH_{\PP^2}) \supseteq
\sigma^*(\rho^* \cL_{\PP^2}).
\]
Indeed, $ \cH_{\PP2}$ is base point free. It is clear that $\rho^*
\cL_{\PP^2}=\cL_{\F_1}+S$, where $S$ is the exceptional curve of
$\rho$. Note that the centers of subsequent blowups  in $\sigma$
(including infinitesimally near centers) lie neither on the proper
transform of $S$ nor on that of general members of $\cL_{\PP^2}$.
Hence,
\[
\sigma^*(\rho^* \cL_{\PP^2})= \sigma^*(\cL_{\F_1}+S)= \cL_{W}+S_W.
\]
Since $S_W$ is $p$-exceptional, applying $p_*$ yields the
assertion. \eproof

\bcor $\cL$ cannot be contained in $|-mK_Y|$ for any $m\in \N$.
\ecor \bproof Assume to the contrary that $\cL\subseteq |-mK_Y|$
for some $m\in \N$. Then the mobile linear system $\cH$, which
defines the birational map $\varphi:Y \dashrightarrow \PP^2$ (see
\ref{lisi}), is also contained in $|-mK_Y|$. This contradicts the
Segre-Manin theorem as stated in \cite[Theorem 2.13]{KSC}. \eproof

The results of this subsection can be summarized as follows.

\bprop\label{sumup} Let $Y\subseteq\PP^3$ be a smooth cubic
surface, and $X\subseteq\A^4$ be the affine cone over $Y$. Then
$X$ admits an effective $\C_+$-action if and only if $Y$ admits a
linear pencil $\cL$ with the following properties: \bnum\item The
base locus $\Bs(\cL)$ consists of a single point, say, $P$, which
is not an Eckardt point on $Y$.
\item A general member $L$ of $\cL$ is singular at $P$,
and $L\setminus\{P\}\simeq\A^1$.
\item $\cL$
has exactly two degenerate members, say $L^{(1)}$ and $L^{(2)}$,
where the curve $L^{(1)}\cup L^{(2)}$ consist of 8 irreducible
components $\Delta_1,\ldots,\Delta_8$.
\item All curves $\Delta_i$, $i=1,\ldots,8$, pass through $P$ and
are pairwise disjoint off $P$. Furthermore,
$\Delta_i\setminus\{P\}\simeq\A^1$ $\forall i$.
\item Every line on $Y$ passing through $P$ is one of the $\Delta_i$.
\item $ -K_Y\equiv D:=\sum_{i=1}^8\delta_i\Delta_i$,
where $\delta_i\in\Q$ and $ 0<\delta_i<1 \,\,\,\forall i$.
\item The pair $(Y,D)$ is not log canonical at $P$.
\item For every $m>0$, $\cL$ is not contained in $|-mK_Y|$.
\enum\eprop

We do not know so far any example of a cubic surface with such a
pencil $\cL$. For the pencils on del Pezzo surfaces from Examples
\ref{misuexa}-\ref{ex11}, not all of the properties (1)-(8) are
fulfilled. For instance, the pencil of Example \ref{ex11}
satisfies (1)-(7), however (8) fails, since $\cL\sim -2K_Y$.

\subsection{The inverse nef value}

\bsit\label{nefbig} The nef value plays an important role in the
adjunction-theoretic classification of polarized varieties. For a
projective variety $Y$ polarized by a nef divisor $H$ we define
the inverse nef value $t_0=t_0(Y,H)$ to be the supremum of $t$
such that the divisor $H+t K_Y$ is nef i.e., \be\label{nefv}
H\cdot C\ge t(-K_Y)\cdot C= t\deg (C)\ee for every curve $C$ on
$Y$. By the Kawamata rationality theorem \cite[Thm.\ 7.1.1]{Mat},
$t_0$ is achieved and is rational. By the Kawamata-Shokurov
base-point-free theorem  \cite[Thm.\ 6.2.1]{Mat}, the divisor
$H+t_0 K_Y$ is semiample i.e., the complete linear system
$|m(H+t_0 K_Y)|$ has no base point for all $m\gg 0$ and defines a
surjective morphism $\varphi\colon Y \rightarrow Y'$ with
connected fibers onto a normal projective variety $Y'$. In
particular, $\kappa(H+t_0K_Y)\ge 0$, where $\kappa$ stands for the
Iitaka-Kodaira dimension.

For a smooth cubic surface $Y$ in $\PP^3$ satisfying Convention
\ref{1000}, we let $\cH=\varphi^{-1}_*(\cO_{\PP^2}(1))$ be the
mobile linear system on $Y$ constructed in \ref{lisi}. In this
setting the inverse nef value $t_0=t_0(Y,\cH)$ is a positive
integer (indeed, for $t=t_0$ the equality in (\ref{nefv}) is
achieved on a $(-1)$-curve). Moreover, $\kappa(\cH+t_0K_Y)=0$ if
and only if $\cH+t_0 K_Y \equiv 0$. However, in the latter case by
Corollary 4.19 a $(-K_Y)$-polar cylinder on $Y$ cannot exist.

By virtue of Theorem \ref{nefna} below, the same conclusion holds
in the case where $\kappa(\cH+t_0K_Y)=1$. In the latter case the
linear system $|m(\cH+t_0K_Y)|$ defines for $m\gg 1$ a conic
bundle $Y\to\PP^1$. Indeed, the image curve is rational since $Y$
is, and an irreducible general fiber $F$ with $F^2=0$ and
$-K_F=-K_Y|_F$ ample is a smooth conic. Actually $|\cH+t_0K_Y|$
defines already a conic bundle. For assuming that
$\cH+t_0K_Y\equiv \beta F$, where $\beta\in\Q$, and taking
intersection with a line $l$ on $Y$ such that $F\cdot l=1$, we
obtain $\beta \in\N$. \esit

\bthm\label{nefna} Let  $\chi:Y\dashrightarrow \PP^2$ be a
birational map and $\cH=\chi_*^{-1} (|\cO_{\PP^2}(1)|)$ be the
proper transform on $Y$ of the complete linear system of lines on
$\PP^2$. Then there is no conic $F$ on $Y$ such that
\be\label{de2} \cH\sim -aK_Y+bF\qquad\mbox{for some}\quad
a\in\N\,\,\,\mbox{and}\,\,\,b\in\Z\,.\ee \ethm

\bproof We use the methods developed in \cite{Is1, Is2}. Consider
a resolution of indeterminacies of $\chi$:

\[
\xymatrix{
& {\widetilde Y}\ar[dl]_{p}\ar[dr]^{q} &&      \\
 Y \ar@{-->}[rr]^{\chi} && \PP^2\\}
\]

\medskip

\noindent Decomposing $p$ into a sequence of blowups with
exceptional curves $E_1,\ldots,E_n$, the linear system $\widetilde
{\cH}=q^*(|\cO_{\PP^2}(1)|)$  on $\widetilde Y$ and the line
bundle $K_{\widetilde Y}$ can be written in $\Pic ({\widetilde
Y})$ as \be\label{de0} \widetilde {\cH}=p^*(\cH)-\sum_{i=1}^n m_i
E_i^*\quad\mbox{and}\quad K_{\widetilde Y} = p^*(K_Y)+\sum_{i=1}^n
E_i^*\,.\ee  Computing the intersection numbers $\widetilde
{\cH}^2$ and $\widetilde {\cH}\cdot K_{\widetilde Y}$, by
(\ref{de0}) we obtain \be\label{de1} 1=\cH^2-\sum_{i=1}^n
m_i^2\quad\mbox{and}\quad -3=K_Y\cdot\cH+\sum_{i=1}^n m_i\,.\ee
Suppose on the contrary that (\ref{de2}) holds for some conic $F$
on $Y$. We choose the minimal possible value of $a>0$. Since
$F^2=0$ on $Y$, from (\ref{de2}) and (\ref{de1}) we deduce
\be\label{de3} \sum_{i=1}^n m_i^2= 3a^2+4ab-1\quad\mbox{and}\quad
\sum_{i=1}^n m_i=3a+2b-3\,.\ee In the rest of the proof we use the
following Claims 1-4.

\medskip

\noindent {\em Claim 1. There is a birational transformation
$Y\dashrightarrow Y'$, where $Y'$ is again a smooth cubic surface
in $\PP^3$, such that for the proper transforms $\cH'$ of $\cH$
and $F'$ of $F$ on $Y'$ we have
$$
\cH'\sim -aK_Y+b'F'$$ with the same $a$ as in (\ref{de2}), and
additionally with
$$m':=\max_i\,\{m_i'\}\le a\quad\forall i\,,$$ where the integers $m_i'$
have the same meaning on $Y'$ as the $m_i$ have on $Y$.}

\medskip

\noindent {\em Proof of Claim 1.} Suppose that $m_i>a$ for some
value of $i$. Consider the conic bundle
$\varphi=\varphi_{|F|}:Y\to\PP^1$. Let us perform an elementary
transformation at the point $P=p(E_i)\in Y$. First we apply the
blowup $\sigma:\widehat Y \to Y$ of $P$ with exceptional divisor
$E$. Assuming that $m_i=\mult_P(\cH)$, on the new surface
$\widehat Y$ we have \be\label{de6} {\widehat
\cH}:=\sigma_*^{-1}(\cH) =\sigma^* (\cH)-m_iE\quad\mbox{and}\quad
K_{\widehat Y} =\sigma^* (K_Y)+E\,.\ee Modulo linear equivalence
we may choose the conic $F$ passing through the point $P$. Then
$F$ is irreducible. Indeed, otherwise $F=F_1+F_2$, where $F_1,F_2$
are two lines on $Y$ and $F_1$ passes through $P$. So
$$a<m_i\le(F_1\cdot\cH)_P\le F_1\cdot\cH=
F_1\cdot(-aK_Y+b(F_1+F_2))=a\,,$$
which is impossible. Thus the proper transform ${\widehat F}\sim
\sigma^*(F)-E$ of $F$  on $\widehat Y$ is a $(-1)$-curve.

The contraction $\sigma':{\widehat Y}\to Y'$ of ${\widehat F}$ to
a point $P'\in Y'$ yields  a smooth conic $F':=\sigma'_*(E)$
passing through $P'$ on the resulting cubic surface $Y'$, such
that
$$\cH':=\sigma'_*(\widehat \cH)\sim -aK_Y+b'F'
\qquad\mbox{for some}\quad b'\in\Z\,.$$
Using (\ref{de2}) and (\ref{de6}), on  $Y'$ we obtain
$$\mult_{P'} (\cH')={\widehat \cH}\cdot {\widehat F}
=\sigma^*(\cH)\cdot {\widehat F}-m_iE\cdot {\widehat F}=\cH\cdot
F-m_i=-aK_Y\cdot F-m_i=2a-m_i<a\,.$$ Iterating this procedure we
achieve finally that $m_i'\le a$ for all values of $i$, as
required. \hfill\qed

So we assume in the sequel that \be\label{de7}
m=\max_i\,\{m_i\}\le a\,\,\forall i\,.\ee

\medskip

\noindent {\em Claim 2. Under the assumption (\ref{de7}) we have
$b<0$.}

\medskip

\noindent {\em Proof of Claim 2.} From (\ref{de3}) and (\ref{de7})
we obtain \be\label{de8} 3a^2+4ab=1+\sum_{i=1}^n m_i^2\le
1+m\sum_{i=1}^n m_i\le 1+m(3a+2b-3)\le 1+a(3a+2b-3)\,.\ee   It
follows by (\ref{de8}) that $2ab\le 1-3a$. Since $a\ge 1$ then
$b\le \frac{1}{2a}-\frac{3}{2}\le -1\,.$ \hfill\qed

\medskip

\noindent {\em Claim 3 (the Noether-Fano Inequality). For $m$ as
in (\ref{de7}) we have $$a\ge m>a+b\,.$$}

\noindent {\em Proof of Claim 3.} The first inequality follows by
(\ref{de7}). To show the second, suppose on the contrary that
\be\label{de4} m\le a+b\,.\ee From (\ref{de8}) and (\ref{de4}) we
obtain \be\label{de121} 3a^2+4ab\le 1+m(3a+2b-3)\le
1+(a+b)(3a+2b-3)\,.\ee Thus by (\ref{de4}) and (\ref{de121})
\be\label{de9} 3\le 3m\le 3(a+b) \le 1+b(a+2b)\,.\ee We claim that
$a+2b\ge 0$. Indeed, let $C$ be the residual line of the conic $F$
on $Y$ so that $F+C\sim -K_Y$. Then by (\ref{de2}),
$$0\le \cH\cdot C=\left(-aK_Y+bF\right)\cdot C=a+2b\,.$$
Now (\ref{de9}) leads to a contradiction, since $b<0$ by Claim 2.
\hfill\qed

\medskip

\noindent {\em Claim 4. Consider the morphism $\varphi:Y\to\PP^1$
defined by the pencil of conics $|F|$ on $Y$. Let $m=m_i$, and let
$Q=p(E_i)\in Y$. If a line $l$ on $Y$ passes though $Q$, then $l$
is a component of the fiber of $\varphi$ through $Q$. }

\medskip

\noindent {\em Proof of Claim 4.} We have \be\label{de11}
(\cH\cdot l)_Q\ge\mult_Q (\cH)= m>a+b\,.\ee On the other hand,
\be\label{de12} (\cH\cdot l)_Q\le\cH\cdot
l=\left(-aK_Y+bF\right)\cdot l=a+bF\cdot l\,.\ee By (\ref{de11})
and (\ref{de12}), $b<bF\cdot l$, where $b<0$ by Claim 2. Therefore
$F\cdot l=0$, and the claim follows. \hfill\qed

\medskip

Let again $C$ be the residual $(-1)$-curve of the conic $F$  on
$Y$ so that $-K_Y\sim F+C$. We let $V$ denote a del Pezzo surface
of degree $4$ obtained by the contraction $\pi : Y\to V$ of $C$,
and $\cH_V,\,F_V$, etc. denote the images on $V$ of $\cH,\,F$,
etc. Due to (\ref{de2}), \be\label{de13} -K_V\sim
F_V\quad\mbox{and}\quad \cH_V\sim -aK_V+bF_V\sim-(a+b)K_V\,.\ee By
Claim 4 there is no line on $V$ through the point $Q_V:=\pi(Q)$.
The blowup $\sigma:Y'\to V$ at $Q_V$ yields yet another cubic
surface $Y'$  with the exceptional $(-1)$-curve
$E=\sigma^{-1}(Q_V)$. For the proper transform
$\cH'=\sigma_*^{-1}(\cH_V)$ on $Y'$ we obtain by (\ref{de13}):
\be\label{de10} \cH'\sim\sigma^*(\cH_V)-mE\sim
(a+b)\sigma^*(-K_V)-mE\sim (a+b)(-K_{Y'})+(a+b-m)E\,.\ee The
linear system
$$|F'|=|-K_{Y'}-E|\,$$ defines a conic
bundle on $Y'$. Plugging $E\sim -K_{Y'}-F'$ into (\ref{de10}) we
deduce: \be\label{de14} \cH'\sim
(2a+2b-m)(-K_{Y'})-(a+b-m)F'\,.\ee Using Claims 2 and 3,
\be\label{de15} 2a+2b-m=(a+b)+(a+b-m)<a\,.\ee  By virtue  of
(\ref{de14}) the latter inequality contradicts the minimality of
$a$. Now the proof of Theorem \ref{nefna} is completed. \eproof

\bcor\label{bigd} Under Convention 4.1 the divisor $\cH+t_0 K_Y$
in \ref{nefbig} is big i.e., $$\kappa(\cH+t_0 K_Y)=2\,.$$ \ecor

\section{Cones over some rational Fano threefolds}
In this section we provide examples of two families of rational
Fano 3-folds such that the affine cones over their anti-canonical
embeddings admit nontrivial $\C_+$-actions.

\bprop\label{quaint} Consider a smooth intersection $Y=Y_{2,2}$ of
two quadric hypersurfaces in $\PP^5$. Then the affine cone $X$
over $Y$ admits an effective $\C_+$-action. \eprop

\bproof According to the criterion of Theorem \ref{ext1}, it is
enough to construct a ($-K_Y$)-polar open cylinder on $Y$. Fix a
line $l\subseteq Y$. Consider the diagram

\[
\xymatrix{
& {\tilde Y}\ar[dl]_{\sigma}\ar[dr]^{\varphi} &&      \\
 Y \ar@{-->}[rr]^{\psi} && \PP^3\\}
\]

\noindent where $\psi$ is the projection with center $l$, $\sigma$
is the blowup of $l$, and $\varphi$ is the blowdown of the divisor
$D$ which is swept out by lines meeting $l$ (i.e., $\sigma(D)$ is
the union of lines meeting $l$; see \cite[Ch. 6]{GH}). It is
easily seen that $\G=\varphi(D)$ is a smooth quintic curve in
$\PP^3$ of genus $2$. The image $Q=\varphi(E)$ of the exceptional
divisor $E$ of $\sigma$ is a quadric  in $\PP^3$. For a line
$l\subseteq Y$, the following alternative holds: either \bnum\item
$N_{l/X}={\mathcal O}_l\oplus {\mathcal O}_l$ and $Q$ is smooth,
\enum or \bnum\item[(2)] $N_{l/X}={\mathcal O}_l(1)\oplus{\mathcal
O}_l(-1)$ and $Q$ is singular.\enum Anyhow,
$$Y\setminus\sigma(D)\simeq\PP^3\setminus Q\,.$$
Suppose that $Q$ is singular; then $Q$ is a quadratic cone. Let
$\Pi$ be a plane in $\PP^3$ passing through the vertex $P$ of $Q$.
We claim that $\PP^3\setminus (Q\cup\Pi)$ is a principal cylinder.
Indeed, consider the projection $\pi_P$ with center $P$ and its
resolution:

\[
\xymatrix{
& {\tilde \PP^3}\ar[dl]_{\sigma'}\ar[dr]^{\varphi'} &&      \\
 \PP^3 \ar@{-->}[rr]^{\pi_P} && \PP^2\\}
\]

\noindent Let $E'\subseteq {\tilde \PP^3}$ be the exceptional
divisor of $\sigma'$, and let $Q'\subseteq {\tilde \PP^3}$ be the
proper transform of $Q$. Then $C=\varphi'(Q')\subseteq\PP^2$ is a
conic, and $E'$ is a section of the $\PP^1$-bundle ${\tilde
\PP^3}\to\PP^2$. Furthermore, $${\tilde
\PP^3}\simeq\PP\left({\mathcal O}_{\PP^2}\oplus{\mathcal
O}_{\PP^2}(1)\right)\,.$$ Letting ${\Pi'}\subseteq {\tilde \PP^3}$
be the proper transform of $\Pi$, the image
$\varphi'(\Pi')=H\subseteq\PP^2$ is a line. We have
$$\PP^3\setminus (Q\cup\Pi)\simeq {\tilde
\PP^3}\setminus (Q'\cup\Pi'\cup E')\,$$ is an $\A^1$-bundle over
$\PP^2\setminus (C\cup H)$. Since $\Pic(\PP^2\setminus (C\cup
H))=0$ we obtain $$\PP^3\setminus (Q\cup\Pi)\simeq\A^1\times
(\PP^2\setminus (C\cup H))\,,$$ as required. \eproof

Let us exhibit yet another family of Fano threefolds with Picard
rank $1$. Their moduli space is $6$-dimensional. Every member $Y$
of this family admits a ($-K_Y$)-polar cylinder, whereas the
subfamily of completions of $\A^3$ is only $4$-dimensional
\cite{Fur}, \cite{Pr1}.

\bprop\label{quaint1} Let $Y=Y_{22}\subseteq\PP^{13}$ be a Fano
variety of genus $12$ with $\Pic (Y)=\Z\cdot[{-K_Y}]$,
anticanonically embedded into $\PP^{13}$. Then the affine cone $X$
over $Y$ admits an effective $\C_+$-action.\eprop

\bproof Again, it is enough to construct a ($-K_Y$)-polar open
cylinder on $Y$. Then the result follows by applying Theorem
\ref{ext1}. Picking a line $l_1\subseteq Y$ we consider the
commutative diagram

\[
\xymatrix{
{\tilde Y}\ar[d]_{\sigma_1}\ar@{-->}[rr]^{\chi} &&
\,\,\,\,\,{\tilde Y}^+\ar[d]^{\varphi_1}     \\
 \,\,\,\,\,\PP^{13}\supseteq Y=Y_{22}\,\,\,\,\, \ar@{-->}[rr]^{\psi_1}
 && \,\,\,\,\,\,\,\,\,\,\,\,\,Y_5\subseteq\PP^6\\}
\]

\medskip

\noindent where  $\sigma_1$ is the blowup of $l_1$, $\psi_1$ is
the double projection with center $l_1$ \footnote{That is,
$\psi_1:Y\dashrightarrow\PP^6$ is defined by the linear system
$|\sigma_1^*\mathcal O_Y(1)-2\tilde E|$, where $\tilde E\subseteq
\tilde Y$ is the exceptional divisor of $\sigma_1$.} onto a Fano
threefold $Y_5$ of degree $5$ and of Fano index $2$,
anticanonically embedded into $\PP^6$, $\varphi_1$ is the blowup
of a smooth rational curve $\G\subseteq Y_5$ of degree $5$, and
$\chi$ is a flop; see \cite[\S 4.3]{IP}. We have
$$Y\setminus H_1\simeq Y_5\setminus H_2\,,$$ where $H_1\subseteq
Y$ is a hyperplane section with $\mult_{l_1}(H_1) =3$, and
$H_2\subseteq Y_5$ is a hyperplane section passing through $\G$.
Thus it suffices to show that $Y_5\setminus H_2$ contains an
$H_2$-polar cylinder.

Let further $l_2\subseteq Y_5$ be a line. Recall that the family
of all lines on $Y_5$ is parameterized by $\PP^2$, and either
$N_{l_2/Y_5}\simeq \mathcal O_{l_2}\oplus\mathcal O_{l_2}$, or
$N_{l_2/Y_5}\simeq \mathcal O_{l_2}(1)\oplus\mathcal O_{l_2}(-1)$.
The lines of second type are parameterized by a smooth conic on
$\PP^2$; see \cite{FN}. There exists a line $l_2$ on $Y_5$ of
second type contained in $H_2$. Consider the projection $\psi_2$
with center $l_2$ and its resolution:

\[
\xymatrix{&
{\tilde Y_5}\ar[dl]_{\sigma_2} \ar[dr]^{\varphi_2}&&     \\
\PP^6\supseteq Y_5\,\,\, \ar@{-->}[rr]^{\psi_2} &&
\,\,\,\,\,Q\subseteq\PP^4\\}
\]

\medskip

\noindent where $\sigma_2$ is the blowup of $l_2$, and
$Q\subseteq\PP^4$ is a smooth quadric. We have $$Y_5\setminus
H_2'\simeq Q\setminus H_3\,,$$ where $H_2'\subseteq Y_5$ and
$H_3\subseteq Q$ are  hyperplane sections such that $\mult_{l_2}
(H_2') =2$, and $H_2'$ is swept out by lines meeting $l_2$. Since
$\psi_2$ is a projection,
$$Y_5\setminus (H_2\cup H_2')\simeq Q\setminus (H_3\cup H_3')\,,$$
where $H_3'\subseteq Q$ is another
 hyperplane section (possibly $H_3'=H_3$).
 It remains to show that the complement
 $Q\setminus (H_3\cup H_3')$ contains an $H_3$-polar cylinder.

 We may assume that
 $H_3'\neq H_3$. The projection $\pi_P:\PP^4\dashrightarrow\PP^3$
 with center at a general point
 $P\in H_3\cap H_3'$ yields an isomorphism
 $$Q\setminus (H_3\cup H_3')\simeq \PP^3\setminus
 (\Pi_1 \cup\Pi_2\cup\Pi_3)\,,$$
 where $\Pi_1, \,\Pi_2,\,\Pi_3$ are three planes in $\PP^3$.
 So the existence
 of an $H_3$-polar cylinder on $Q\setminus (H_3\cup H_3')$
 is equivalent to the existence
 of a $\Pi_1$-polar cylinder on $\PP^3\setminus
 (\Pi_1 \cup\Pi_2\cup\Pi_3)$.
 Now the assertion easily follows.
\eproof

\end{document}